\documentclass[moor,nonblindrev]{informs1}

\usepackage{physics}
\usepackage{mathtools}
\usepackage{tensor}

\newcommand{\ceil}{\mathop{\rm ceil} \nolimits}
\newcommand{\floor}{\mathop{\rm floor} \nolimits}
\newcommand{\half}{\tfrac{1}{2}}
\newcommand{\asymx}{\mathop{\sim}}
\newcommand{\asym}[1]{\mathrel{\asymx_{#1}}}
\newcommand{\diag}{\mathop{\rm diag} \nolimits}
\newcommand{\sgn}{\mathop{\rm sgn} \nolimits}

\newcommand{\erfc}{\mathop{\rm erfc} \nolimits}
\newcommand{\erfcx}{\mathop{\rm erfcx} \nolimits}

\newmuskip\pFqmuskip
\newcommand*\pFq[6][8]{%
  \begingroup 
  \pFqmuskip=#1mu\relax
  \mathcode`\,=\string"8000
  \begingroup\lccode`\~=`\
  \lowercase{\endgroup\let~}\pFqcomma
  {}_{#2}F_{#3}{\left[\genfrac..{0pt}{}{#4}{#5};#6\right]}%
  \endgroup
}
\newcommand{\pFqcomma}{\mskip\pFqmuskip}

\def\wscl{0.8}
\def\hscl{0.56}

\usepackage{natbib}
\NatBibNumeric
\def\bibfont{\small}%
\bibpunct[, ]{[}{]}{,}{n}{}{,}%

\definecolor{myblue}{RGB}{0, 20, 114}
\usepackage[colorlinks=true,breaklinks=true,bookmarks=true,urlcolor=myblue,
     citecolor=myblue,linkcolor=myblue,bookmarksopen=false,draft=false]{hyperref}

\def\EMAIL#1{\href{mailto:#1}{#1}}

\TheoremsNumberedThrough     

\EquationsNumberedThrough    

\MANUSCRIPTNO{0} 

\def\theARTICLETOP{}
\def\theLRHSecondLine{}
\def\theRRHSecondLine{}

\begin{document}


\RUNAUTHOR{Zuk and Kirszenblat}

\RUNTITLE{Partial Busy Period}

\TITLE{Exact Results for the Distribution of the Partial Busy Period for a Multi-Server Queue}

\ARTICLEAUTHORS{%
\AUTHOR{Josef Zuk}
\AFF{Defence Science and Technology Group, Melbourne, Australia,
\EMAIL{josef.zuk@defence.gov.au}}
\AUTHOR{David Kirszenblat}
\AFF{Defence Science and Technology Group, Melbourne, Australia,
\EMAIL{david.kirszenblat@defence.gov.au}}
} 

\ABSTRACT{%
Exact explicit results are derived for the distribution of the partial busy period of the M/M/$c$ multi-server
queue for a general number of servers.
A rudimentary spectral method leads to a representation that is amenable to efficient numerical computation
across the entire ergodic region.
An alternative algebraic approach yields a representation as a finite sum of Marcum Q-functions depending
on the roots of certain polynomials that are explicitly determined for an arbitrary number of servers.
Asymptotic forms are derived in the limit of a large number of servers under two scaling regimes,
and also for the large-time limit. Connections are made with previous work.
The present work is the first to offer tangible exact results for the distribution when the number of servers
is greater than two.
}%

\MSCCLASS{Primary: 90B22; secondary: 60K25, 60J74}
\ORMSCLASS{Primary: Queues: Busy period analysis; secondary: Queues: Algorithms}

\HISTORY{Date created: May 27, 2023. Last update: August 21, 2023.}

\maketitle

%

\section{Introduction}
One of the fundamental concepts in queueing theory is the busy period (BP).
Other basic quantities include the stationary queue-size and waiting-time distributions.
For the archetypal M/M/$c$ multi-server queueing model,
the latter are extensively covered in all introductory textbooks,
{\it e.g.}\ \citep{BP:Gnedenko89, BP:Shortle18}.
But they are relevant only in the ergodic region where the total traffic intensity is
less than unity. By contrast,
the BP persists as a valid concept beyond this region.
However, there is little available to be said about its calculation,
even in the ergodic case.

There are two main characterizations of a BP.
The full BP refers to the time interval during which all servers are occupied.
The partial BP refers to the time interval during which at least one
server is occupied.
For a single-server queue, the two definitions trivially coincide.
It is the distribution of the partial BP that will be studied here, and use of the
term `busy period' on its own will assume the partial variant.

We shall confine our attention to the M/M/$c$ queue,
and denote the number of servers by $N$. Thus
\mbox{$c = N$}
in our discussion.
Within this class of models, the BP distribution for the single-server case is known.
For a multi-server queue, the full BP is simply related to that of the
single-server case with an adjusted service time \citep{BP:Daley98}.
On the other hand, results on the partial BP for a multi-server queue are sparse.

For more than two servers, there exists in the published literature neither
numerical algorithms nor analytic expressions for the distribution of the
partial BP. The present work provides both of these.
We first demonstrate a spectral method that serves as the basis for a
simple and efficient numerical algorithm that covers at least the entire ergodic region
and handles server numbers up to around $N = 80$. We have tested it on the interval
of traffic intensities
\mbox{$0 \leq r \leq 1$}.
We also develop an algebraic method that gives rise to an explicit representation of
the BP distribution for any number of servers
as a finite sum of Marcum Q-functions, dependent on the
roots of a certain family of polynomials whose coefficients we determine.

\citet{BP:Karlin58} give an explicit integral representation for the case of two servers
based on a spectral method involving families of orthogonal polynomials.
But, while they claim that their method can in principle be extended to a larger number of servers,
this must be done on a case by case basis, and the calculations
become so unwieldy beyond
\mbox{$N = 2$}
that this line of attack has not hitherto been pursued by anyone.
By contrast, we show that a much more rudimentary spectral method yields
more generally applicable results with greater ease.
\citet{BP:Arora64} has also derived the distribution for the two-server problem
as a infinite series of modified Bessel functions.
Our algebraic approach can be viewed as a generalization of this to an arbitrary
number of servers.

All other work appears to be centered around computing moments of the BP distribution.
\citet{BP:Natvig75} obtains the first and second order moments of the length of the partial BP.
Second moments of the partial BP distribution are also derived by \citet{BP:Omahen78}.
Further progress has been made more recently by \citet{BP:Artalejo01} who show how to compute arbitrary moments.
As commented upon by these authors, prior to their paper in 2001, general results were only known for the one
and two server cases, and moments up to second order had appeared for the general problem.
We are not aware of anything more recent that alters this state of affairs.

\section{Motivation}
The BP distribution for the basic M/M/$c$ queue is equally applicable to a large family
of more complex, albeit memoryless,
queueing systems that incorporate multiple classes of arrivals each comprised of
multiple priority levels that are processed under a variety of queuing disciplines
({\it e.g.}\ preemptive or non-preemptive).
A concrete example is furnished by the ambulance ramping problem encountered by
hospital emergency departments (ED), where one question that may be studied is to what extent an
ambulance offload zone alleviates the ambulance queue that builds up while patients wait to be
admitted and prevents ambulances from being dispatched to other calls \citep{BP:Li19,BP:Hou20}.
Patients are differentiated via arrival source into the ED as either walk-ins and ambulance arrivals,
with ambulance arrivals further sub-divided into those who enter the ED from
the ambulance queue or the offload zone.
Multiple priorities are assigned to patients in each arrival class depending on their acuity levels,
and each arrival-class/priority-level combination can exhibit a different arrival rate.
The main constraint is that of a common mean service time
({\it i.e.}\ treatment time) for all patients.
It does not matter to the servers how the patients are arranged prior to entry,
or on the mixture currently in service.
Hence the BP distribution is just that of the basic model.
\citet{BP:Rastpour15} has recently studied partial BP distributions in the context of
emergency medical services.

Our interest in the partial BP emerged from its
relevance to regenerative simulation of the steady-state limit
of multi-server, multi-priority, multi-class queues as discussed above.
The regenerative simulation technique was introduced by \citet{BP:Crane74}.
Empirical distributions for queue lengths and waiting times
per class or priority level can be ascertained
from steady-state discrete event simulations, and subsequently compared with
theoretical predictions.
In the ergodic region, sample means for summary statistics are easily estimated.
However, confidence intervals are also required to judge whether
the null hypothesis that the empirical and theoretical values are generated by the
same candidate distribution is to be rejected at some given level of statistical
significance.
In steady state simulation, since one is producing a very long single run of a
stochastic ({\it e.g.}\ Markov) process, the data underlying the empirical distributions
are highly correlated.
Difficulties with estimating confidence intervals arising from unwanted correlation
can be avoided by recognizing that one is dealing with a regenerative process,
and partitioning the time series into consecutive regeneration cycles. These are
statistically independent and identically distributed as each cycle effectively restarts
the process without memory of the past. The natural regeneration point in a queuing problem
is the empty state, which occurs an infinite number of times in an ergodic system.
Thus, each regeneration cycle comprises an initial idle ({\it i.e.}\ empty) period
followed by a partial BP.
The distribution of the idle period is simply the inter-arrival distribution,
assumed here to be exponential.
Therefore, the interesting quantity in the regeneration cycle problem is the
partial BP.
In carrying out a state-state simulation of a queueing system, it is a straightforward
matter to collect data on the lengths of regeneration cycles and BPs.
Comparison of these empirical results with theoretical expectations
serves as a useful diagnostic tool for judging the soundness of the simulation.
It was the paucity of theoretical results for this comparison that led to the present investigation.

The remainder of the paper is organized as follows.
Our starting point in Section~\ref{Spectral} will be the recursive system of equations
considered by \citet{BP:Artalejo01},
which can be solved to give the moment generating function (MGF) of the partial BP.
We begin by recasting this as a matrix inversion problem for a tridiagonal matrix.
The matrix, which has the form of a generator for a birth-death process, is inverted
by means of its spectral decomposition.
This leads to an explicit representation of the distribution in terms of the eigenstates
that is amenable to efficient numerical implementation.
This is followed in Section~\ref{Algebraic} by a algebraic approach where the
recurrence relations of \citet{BP:Artalejo01}
are solved directly for the MGF. The method leads to a finite continued fraction representation similar
to previous findings in the literature, but which does not appear to be suitable for further analysis.
An alternative method yields an explicit form for the MGF in terms of a family of polynomials, whose degree
increases with the number of servers, and which depend parametrically on the traffic intensity.
A simple two-dimensional system of recurrence relations
allows the coefficients of these polynomials to be computed.
In Section~\ref{Asym}, asymptotic limits as the number of servers $N$ approaches infinity are derived
under two different scaling regimes. These are subsequently used for comparison in validating
exact results for large but finite $N$, and determining how large $N$ must be for one to judge that
asymptotic behaviour has effectively set in.
In Section~\ref{Complex}, the polynomial-based representation of the MGF is combined with a
complex contour integral that implements the inverse Laplace transform,
whose structure is determined by poles that arise from zeros of the polynomials.
The contour integration evaluates
to explicit closed-form expressions for the partial BP distribution
as a finite sum of Marcum Q-functions.
In Section~\ref{Statistics},
various summary statistics are computed from the posited theoretical distributions
and compared, across a large range of model parameters,
with known exact results and with the output of Monte Carlo simulation.
Excellent agreement is demonstrated over the entire ergodic range of traffic intensity,
and for server numbers up to well within the asymptotic regime.
Concluding remarks are presented in Section~\ref{Conclusions}.

\section{Spectral Method}
\label{Spectral}
We consider the M/M/$c$ queue with
\mbox{$c = N$}
servers, Poisson arrival rate $\lambda$ and mean treatment (or service) time $1/\mu$.
Following \citep{BP:Artalejo01}, the quantities
\mbox{$\phi_n(s)$},
\mbox{$n = 1,2,\ldots$},
where $n$ labels the number of busy servers,
are defined to be MGFs of the first passage times
from the state
\mbox{$n = 1,2,\ldots,N$}
to the empty state $0$.
Also, let $T_{\text{bp}}$ be the random variable (RV) describing the partial BP,
in which case its MGF is given by
\mbox{$\phi_1(s) \equiv \langle e^{-sT_{\text{bp}}}\rangle_{T_{\text{bp}}}$}.
Then, if we set
\begin{align}
\begin{aligned}
\psi_\pm(s) &\equiv \frac{1}{2\lambda}\left[\lambda + N\mu + s \pm\sqrt{(\lambda + N\mu + s)^2 - 4N\mu\lambda}\right] \\
     &\equiv \frac{1}{2r}\left[r+ 1 + s/(N\mu) \pm\sqrt{(r + 1 + s/(N\mu))^2 - 4r}\right] \;,
\end{aligned}
\end{align}
where
\mbox{$r \equiv \lambda/(N\mu)$}
is the total traffic intensity,
the MGF for the single-server problem
(\mbox{$N = 1$})
reads
\mbox{$\phi_1(s) = \psi_-(s)$}.
We confine our attention to the ergodic region
\mbox{$0 \leq r < 1$},
but we shall also study the boundary value
\mbox{$r = 1$}
as a special case.
We also introduce
\mbox{$\mu_n \equiv n\mu$},
for
\mbox{$n = 1,2,\ldots,N$},
and without loss of generality, we may set
\mbox{$\mu = 1/N$},
so that
\mbox{$r = \lambda$}
and
\begin{equation}
\psi_\pm(s) = \frac{1}{2r}\left[r + 1 + s \pm\sqrt{(r + 1 + s)^2 - 4r}\right] \;.
\label{PsiNorm}
\end{equation}
Finally, we write
\mbox{$\mu_n \equiv n\mu_1$}
with
\mbox{$\mu_1 = 1/N$},
so that
\mbox{$0 \leq \mu_n \leq 1$}
for all $n$.

The linear recurrence relations of \citet{BP:Artalejo01} may be cast as follows:
For
\mbox{$n = 1,2,\ldots,N-1$},
\begin{equation}
\mu_n\phi_{n-1} - (r + \mu_n + s)\phi_n + r\phi_{n+1} = 0 \;,
\label{PhiRecur}
\end{equation}
and
\mbox{$\phi_N = \psi_-(s)\phi_{N-1}$},
with
\mbox{$\phi_0 \equiv 1$}.
As the original paper does not make the derivation of this result explicit,
we supply the details in the Appendix.
These equations have an $(N-1)\times(N-1)$ tridiagonal matrix structure that may be expressed as
\begin{equation}
\begin{bmatrix}
  r + \mu_1 + s & -r            & 0             & \cdots      & \cdots             & 0                  \vphantom{\phi_1}     \\
  -\mu_2        & r + \mu_2 + s & -r            & \cdots      & \cdots             & 0                  \vphantom{\phi_2}     \\
  0             & -\mu_3        & r + \mu_3 + s & -r          & \cdots             & 0                  \vphantom{\vdots}     \\
  \vdots        &               & \ddots        &  \ddots     &  \ddots            & \vdots             \vphantom{\vdots}    \\
  0             & \cdots        & \cdots        &  -\mu_{N-2} & r + \mu_{N-2} + s  & -r                 \vphantom{\phi_{N-2}} \\
  0             & \cdots        & \cdots        &  0          & -\mu_{N-1}         &  r + \mu_{N-1} + s \vphantom{\phi_{N-1}}
\end{bmatrix} \!\!
\begin{bmatrix}
  \phi_1 \\
  \phi_2 \\
  \vdots \\
  \vdots \\
  \phi_{N-2} \\
  \phi_{N-1}
\end{bmatrix} =
\begin{bmatrix}
  \mu_1 \\
  0 \\
  \vdots \\
  \vdots \\
  0 \\
  r\phi_{N}
\end{bmatrix} \;.
\label{MatEq}
\end{equation}
In Dirac notation \citep{BP:Wiki}, this matrix equation reads
\begin{equation}
(s\mathbb{I} - M)|\phi\rangle = |w\rangle \;,
\label{Dirac}
\end{equation}
where the matrix $M$ is the birth-death process generator
\begin{equation}
M =\begin{bmatrix}
  -(r + \mu_1) & r            & 0            & \cdots     & \cdots            & 0 \\
  \mu_2        & -(r + \mu_2) & r            & \cdots     & \cdots            & 0 \\
  0            & \mu_3        & -(r + \mu_3) & r          & \cdots            & 0 \\
  \vdots       &              & \ddots       &  \ddots    &  \ddots           & \vdots \\
  0            & \cdots       & \cdots       &  \mu_{N-2} & -(r + \mu_{N-2})  & r \\
  0            & \cdots       & \cdots       &  0         & \mu_{N-1}         & -( r + \mu_{N-1})
\end{bmatrix} \;.
\label{TriDiagM}
\end{equation}
It may be compared with the matrices introduced in (1.3) and (4.1) of \citet{BP:Karlin58}.
Its trace and determinant, useful for checking the numerical solution of the eigenvalue
problem, are
\begin{equation}
\tr(M) = -(N-1)\left(r + \half\right) \;, \quad
\det(M) = (-r)^{N-1}\sum_{k=0}^{N-1}\frac{k!}{(Nr)^k} \;.
\end{equation}
If we introduce the standard orthonormal Cartesian basis on $\mathbb{R}^{N-1}$ as
\mbox{$|e_i\rangle$},
\mbox{$i = 1,2,\ldots,N-1$},
then we have that
\mbox{$\phi_n \equiv \langle e_n|\phi\rangle$}.
Also,
\begin{equation}
|w\rangle \equiv \mu_1|e_1\rangle + r\phi_N|e_{N-1}\rangle \;.
\end{equation}

Since $M$ is non-symmetric, it has distinct right and left eigenvectors so that its
eigenvalue problem is written as
\begin{equation}
M|v^{(k)}_{\text{R}}\rangle  = \xi_k|v^{(k)}_{\text{R}}\rangle \;, \quad
     \langle v^{(k)}_{\text{L}}|M = \xi_k\langle v^{(k)}_{\text{L}}| \;.
\end{equation}
With the eigenvectors required to satisfy the orthonormality conditions
\mbox{$\langle v^{(k)}_{\text{L}}|v^{(\ell)}_{\text{R}}\rangle =
     \langle v^{(k)}_{\text{R}}|v^{(\ell)}_{\text{L}}\rangle = \delta_{k\ell}$},
the spectral decomposition of $M$ reads
\begin{equation}
M = \sum_{k=1}^{N-1} |v^{(k)}_{\text{R}}\rangle \xi_k \langle v^{(k)}_{\text{L}}| \;.
\end{equation}
To determine $\phi_1(s)$, we invert (\ref{Dirac})
to obtain
\begin{equation}
|\phi\rangle = \left(s\mathbb{I} - M\right)^{-1}|w\rangle \;,
\label{BraKet}
\end{equation}
and observe that
it suffices to consider just the first and last elements of this equation, namely
\begin{align}
\begin{aligned}
\langle e_1|\phi\rangle &= \mu_1\sum_{k=1}^{N-1}\frac{1}{s - \xi_k}
     \langle e_1|v^{(k)}_{\text{R}}\rangle\langle v^{(k)}_{\text{L}}|e_1\rangle
     + r\phi_N\sum_{k=1}^{N-1}\frac{1}{s - \xi_k}
          \langle e_1|v^{(k)}_{\text{R}}\rangle\langle v^{(k)}_{\text{L}}|e_{N-1}\rangle \;, \\
\langle e_{N-1}|\phi\rangle &= \mu_1\sum_{k=1}^{N-1}\frac{1}{s - \xi_k}
     \langle e_{N-1}|v^{(k)}_{\text{R}}\rangle\langle v^{(k)}_{\text{L}}|e_1\rangle
     + r\phi_N\sum_{k=1}^{N-1}\frac{1}{s - \xi_k}
          \langle e_{N-1}|v^{(k)}_{\text{R}}\rangle\langle v^{(k)}_{\text{L}}|e_{N-1}\rangle \;.
\end{aligned}
\label{BraKetElem}
\end{align}

\begin{figure}
\FIGURE
{\includegraphics[width=\wscl\linewidth, height=\hscl\linewidth]{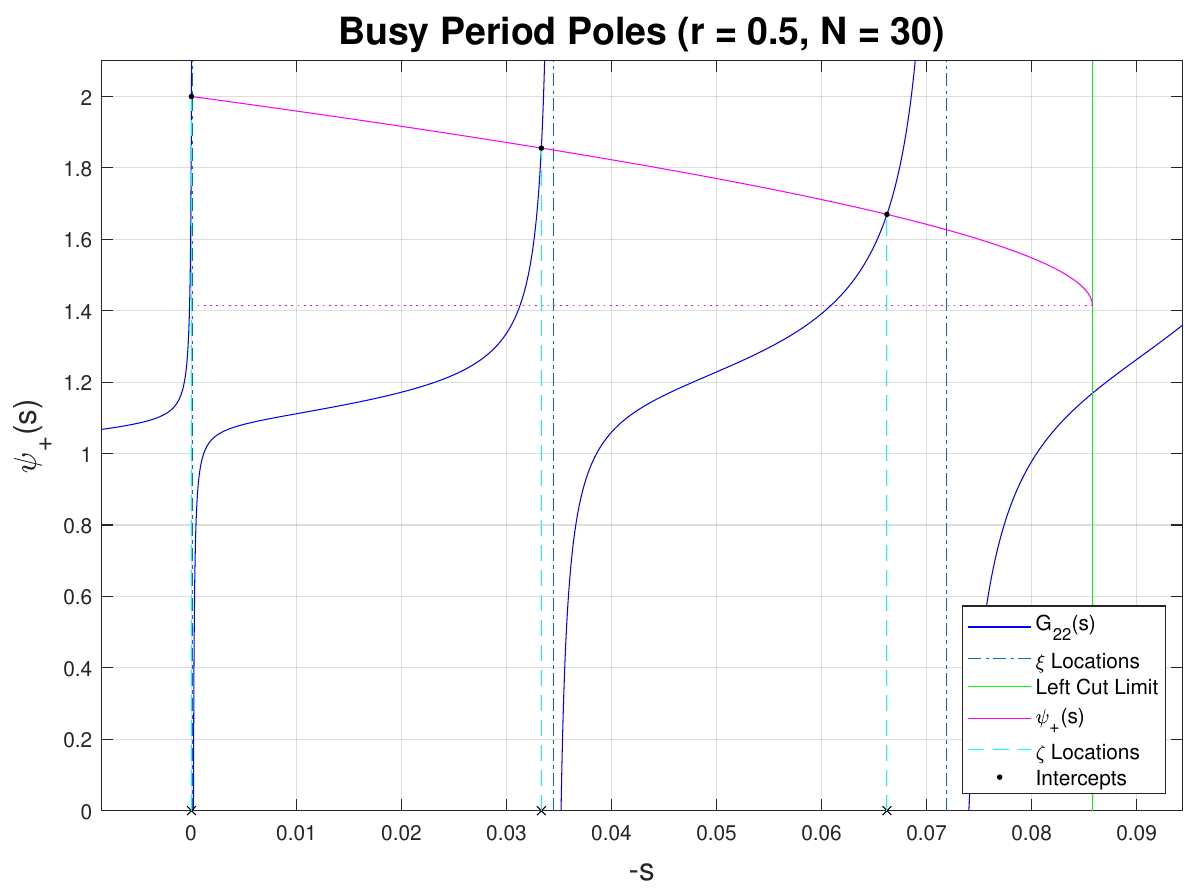}}
{\hphantom{x}\label{BPPoles}}
{BP pole locations for case $r = 0.5$, $N = 30$ at $s = \zeta_\ell$,
     $\ell = 1,2,3$ are indicated by the crosses on the horizontal axis.
     The smallest pole lies very close to zero.}
\end{figure}

The matrix $M$ may be symmetrized by means of a diagonal similarity transformation
\mbox{$D = \diag[d_1,d_2,\ldots,d_{N-1}]$}
where
\mbox{$d_1 = 1$}
and
\mbox{$d_k = \sqrt{\mu_k/r}{\cdot}d_{k-1}$},
for
\mbox{$k = 2,\ldots,N-1$},
according to
\mbox{$S = D^{-1}MD$}.
Its spectral decomposition reads
\begin{equation}
S = \sum_{k=1}^{N-1} |v^{(k)}\rangle \xi_k \langle v^{(k)}| \;,
\end{equation}
and the eigenvectors (assumed orthonormal) are related by
\begin{equation}
|v^{(k)}_{\text{R}}\rangle = D|v^{(k)}\rangle \;, \quad
|v^{(k)}_{\text{L}}\rangle = D^{-1}|v^{(k)}\rangle \;.
\end{equation}
Let us define
\begin{align}
\begin{aligned}
v^{(k)}_{11} &\equiv \langle e_1|v^{(k)}\rangle\langle v^{(k)}|e_1\rangle \;, \\
v^{(k)}_{22} &\equiv \langle e_{N-1}|v^{(k)}\rangle\langle v^{(k)}|e_{N-1}\rangle \;, \\
v^{(k)}_{12} &\equiv \langle e_1|v^{(k)}\rangle\langle v^{(k)}|e_{N-1}\rangle \;, \\
             &=      \langle e_{N-1}|v^{(k)}\rangle\langle v^{(k)}|e_1\rangle \;.
\end{aligned}
\end{align}
It follows that
\begin{align}
\begin{aligned}
\langle e_1|v^{(k)}_{\text{R}}\rangle\langle v^{(k)}_{\text{L}}|e_1\rangle          &= v^{(k)}_{11} \;, \\
\langle e_{N-1}|v^{(k)}_{\text{R}}\rangle\langle v^{(k)}_{\text{L}}|e_{N-1}\rangle  &= v^{(k)}_{22} \;, \\
\langle e_1|v^{(k)}_{\text{R}}\rangle\langle v^{(k)}_{\text{L}}|e_{N-1}\rangle      &= \Lambda^{1/2} v^{(k)}_{12} \;, \\
\langle e_{N-1}|v^{(k)}_{\text{R}}\rangle\langle v^{(k)}_{\text{L}}|e_{1}\rangle    &= \Lambda^{-1/2} v^{(k)}_{12} \;,
\end{aligned}
\end{align}
where
\begin{equation}
\Lambda \equiv \left(d_{1}/d_{N-1}\right)^2 = \prod_{k=2}^{N-1}(r/\mu_k) \;.
\label{Lambda}
\end{equation}
It is useful to note that the following consequence of the completeness relation for the eigenvector basis
\begin{equation}
\sum_{k=1}^{N-1} |v^{(k)}\rangle\langle v^{(k)}| = \mathbb{I} \quad \Rightarrow \quad
\sum_{k=1}^{N-1} v^{(k)}_{mn} = \delta_{mn} \;,
\end{equation}
for
\mbox{$m,n = 1,2$},
where $\delta_{mn}$ is the Kronecker delta.
We can now express (\ref{BraKetElem}) as
\begin{align}
\begin{aligned}
\phi_1     &= \mu_1G_{11} + r\phi_N\Lambda^{1/2}G_{12} \;, \\
\phi_{N-1} &= \mu_1\Lambda^{-1/2}G_{12} + r\phi_N G_{22} \;,
\end{aligned}
\label{PhiPair}
\end{align}
where we have introduced the resolvent functions
\begin{equation}
G_{11}(s) \equiv \sum_{k=1}^{N-1}\frac{v^{(k)}_{11}}{s - \xi_k} \;, \quad
G_{12}(s) \equiv \sum_{k=1}^{N-1}\frac{v^{(k)}_{12}}{s - \xi_k} \;, \quad
G_{22}(s) \equiv \sum_{k=1}^{N-1}\frac{v^{(k)}_{22}}{s - \xi_k} \;.
\end{equation}
We can eliminate $\phi_N$ by applying the relation
\mbox{$\phi_N = \psi_-(s)\phi_{N-1}$}
to the second equation in (\ref{PhiPair}), which yields
\begin{equation}
\phi_N = \frac{\mu_1}{r\Lambda^{1/2}}{\cdot}\frac{G_{12}}{\psi_+ - G_{22}} \;,
\end{equation}
upon noting that
\mbox{$\psi_+(s)\psi_-(s) = 1/r$}.
This expression is then inserted into the first equation in (\ref{PhiPair})
to provide the explicit result for the MGF
\begin{equation}
\phi_1(s)/\mu_1 = G_{11}(s) + \frac{G_{12}^2(s)}{\psi_+(s) - G_{22}(s)} \;.
\label{PhiG}
\end{equation}
It should be noted that, due to cancellations, there are actually no singularities
at the eigenvalues
\mbox{$s = \xi_k$}.
The singularity structure in the complex $s$-plane of $\phi_1(s)$ comprises
a cut and potentially a finite number of poles.
The cut is due to the presence of $\psi_+(s)$, and lies
along the negative real axis corresponding to the finite interval for
which the discriminant
\mbox{$\Delta(s) \equiv b^2(s) - 4r$}
is negative, and where
\mbox{$b(s) \equiv r + 1 +s$}.
Thus, the cut spans
\mbox{$x_- \leq -s \leq x_+$},
where the cut limits are given by
\mbox{$x_\pm = (1 \pm \sqrt{r})^2$}.
The poles arise from the vanishing of the denominator in (\ref{PhiG}), and only lie
on the negative real axis between the origin and the near cut boundary,
as discussed in \citep{BP:Karlin58}.
This is illustrated in Figure~\ref{BPPoles} for the case
\mbox{$r = 0.5$},
\mbox{$N = 30$},
the details of which will be explained later on.

The BP distribution is recovered from the MGF $\phi_1(s)$ by an inverse Laplace transform
implemented as a Bromwich contour integral in the complex $s$-plane.
The contour may be deformed onto the negative real axis, giving rise to
discrete residue contributions from any relevant poles, plus a real-valued integral on the cut.
An explicit representation of this integral has been given for the two-server case in equation (6.25)
of \citep{BP:Karlin58}. It is also of the same type as the integral arising in the waiting-time
distribution for priority queues, as can be seen in \citep{BP:Davis66}.
One should observe that, in the present formalism, the two-server problem is trivial as
the matrix $M$ is one dimensional, so that there is no eigenvalue problem to solve.

To evaluate the contribution to the BP distribution from the cut,
we introduce the product function
\begin{equation}
\Pi(s) \equiv \prod_{k=1}^{N-1} (s - \xi_k) \;,
\end{equation}
whose square will multiply both the numerator and denominator in (\ref{PhiG})
in order to eliminate spurious singularities at the eigenvalues.
This results in the appearance of the non-singular functions
\begin{align}
\begin{aligned}
H_{12}(s) &\equiv G_{12}(s)\Pi(s) &= \sum_{k=1}^{N-1}v^{(k)}_{12}\prod_{\ell=1 \atop \ell\neq k}^{N-1} (s - \xi_\ell) \;, \\
H_{22}(s) &\equiv G_{22}(s)\Pi(s) &= \sum_{k=1}^{N-1}v^{(k)}_{22}\prod_{\ell=1 \atop \ell\neq k}^{N-1} (s - \xi_\ell) \;.
\end{aligned}
\end{align}
Next, letting
\mbox{$D^{\pm}(s) \equiv \psi_\pm(s) - G_{22}(s)$}
so that
\mbox{$D^+(s)$}
is the denominator in (\ref{PhiG}),
we see that
\begin{equation}
D^+(s)D^-(s) = \left[1 - b(s)G_{22}(s) +rG^2_{22}(s)\right]/r \;.
\end{equation}
By appealing to the identity
\mbox{$1/D^+(s) = D^-(s)/[D^+(s)D^-(s)]$},
we are able to construct the cut contribution as
\begin{equation}
\phi_1(s)/\mu_1 \underset{\text{cut}}{\longrightarrow}
     (\pm i)R_{\text{cut}}(s)\sqrt{\left|\Delta(s)\right|} \;,
\end{equation}
where the cut function
\mbox{$R_{\text{cut}}(s)$}
is given by
\begin{equation}
R_{\text{cut}}(s) \equiv \frac{H_{12}^2(s)}{\Pi^2(s) - b(s)\Pi(s)H_{22}(s) + rH_{22}^2(s)} \;.
\end{equation}
This definition takes account of the fact that the cut is traversed in two opposite
directions.
For
\mbox{$N = 1$},
the cut function is just the constant
\mbox{$R_{\text{cut}}(s) = 1/r$}
while, for
\mbox{$N = 2$},
we obtain
\begin{equation}
R_{\text{cut}}(s) = \frac{2}{r - 1/2 - s} \;.
\end{equation}

Now, since
\mbox{$\Pi(\xi_k) = 0$}
for all
\mbox{$k = 1,2,\ldots,N-1$},
it follows that
\begin{equation}
R_{\text{cut}}(\xi_k) = \frac{1}{r}{\cdot}\left[\frac{H_{12}(\xi_k)}{H_{22}(\xi_k)}\right]^2 \;.
\end{equation}
But we also have that
\begin{equation}
H_{12}(\xi_k) = v^{(k)}_{12}\prod_{\ell=1 \atop \ell\neq k}^{N-1} (\xi_k - \xi_\ell)\;, \quad
H_{22}(\xi_k) = v^{(k)}_{22}\prod_{\ell=1 \atop \ell\neq k}^{N-1} (\xi_k - \xi_\ell)\;.
\end{equation}
Thus,
\begin{equation}
R_{\text{cut}}(\xi_k) = (v^{(k)}_{12}/v^{(k)}_{22})^2/r \;,
\end{equation}
which provides a useful numerical check.
By considering the point
\mbox{$s = 0$},
we obtain the identities
\begin{align}
\begin{aligned}
1 &= \mu_1G_{11}(0) + rG_{12}\sqrt{\Lambda} \;, \\
1 &= \mu_1G_{11}(0) + \mu_1G_{12}^2(0)/[1/r - G_{22}(0)] \;,
\end{aligned}
\end{align}
which imply that
\begin{equation}
\mu_1G_{12}(0)/\sqrt{\Lambda} + rG_{22}(0) = 1 \;.
\label{G12G22}
\end{equation}
The eigenvalue problem for the symmetric matrix $S$ is solved numerically using the LAPACK
routine {\sf dstevd}, that implements a divide-and-conquer method \citep{BP:LAPACK99}.
Identical results are produced when using {\sc Matlab}'s {\sf eig}
function, which is probably an indication that {\sf eig} invokes the same algorithm when
applied to this problem.

On making the change of integration variable
\mbox{$x = -s$}
to obtain positive values for the coordinates of the cut, we can express the cut contribution to the BP
probability density function (PDF) as
\begin{equation}
P_{\text{cut}}(t) = \frac{\mu_1}{2\pi}\int_{(1-\sqrt{r})^2}^{(1+\sqrt{r})^2} dx\,
     e^{-xt}R_{\text{cut}}(-x)\sqrt{|\Delta(-x)|} \;.
\end{equation}
For the survival function\footnote{Also known as the complementary cumulative distribution function.} (SF),
\begin{equation}
\bar{F}_{\text{cut}}(t) = \frac{\mu_1}{2\pi}\int_{(1-\sqrt{r})^2}^{(1+\sqrt{r})^2} \frac{dx}{x}\,
     e^{-xt}R_{\text{cut}}(-x)\sqrt{|\Delta(-x)|} \;.
\end{equation}
Let
\mbox{$\Delta x \equiv x_+ - x_- = 4\sqrt{r}$}
denote the length of the cut. Then a further change of variable, such that
\mbox{$ x = x_- + \Delta x{\cdot}u$},
puts the integral over the unit interval
\begin{equation}
\bar{F}_{\text{cut}}(t) = \frac{2\mu_1\sqrt{r}}{\pi}e^{-4a\sqrt{r}t}\int_0^1
     \frac{du}{u + a}\, e^{-4\sqrt{r}tu}R_{\text{cut}}(-4\sqrt{r}(u+a))\sqrt{u(1-u)} \;,
\label{FcutU}
\end{equation}
where
\begin{equation}
a \equiv x_-/\Delta x = \tfrac{1}{4}\left(\sqrt{r} + 1/\sqrt{r}\right) - \tfrac{1}{2} \;.
\end{equation}
A further change of variable
\mbox{$v = 2u - 1$}
casts the integral into a form that is amenable to efficient Gauss-Chebyshev quadrature:
\begin{equation}
\bar{F}_{\text{cut}}(t) = \frac{\mu_1\sqrt{r}}{\pi}e^{-(r+1)t}\int_{-1}^{+1} dv\, \sqrt{1 - v^2}
     e^{-2\sqrt{r}tv}{\cdot}\frac{R_{\text{cut}}(-2\sqrt{r}(v+\alpha))}{v+\alpha} \;,
\label{FcutV}
\end{equation}
with
\mbox{$\alpha \equiv 2a+1 = (\sqrt{r} + 1/\sqrt{r})/2$}.
It is known that Gauss-Chebyshev quadrature of the second kind
is equivalent to the trapezoidal rule on the unit circle \citep{BP:Chawla70}.
Thus, if we set
\mbox{$v = \cos(\pi\tau)$},
so that
\mbox{$0 \leq \tau \leq 1$},
and subsequently uniformly discretize the unit interval into $L$ sub-intervals according to
\mbox{$\tau_n = n/L$},
then we obtain the $L$-point quadrature rule for the SF given by
\begin{equation}
\bar{F}_{\text{cut}}(t) \simeq \sum_{n=1}^L w_n e^{-x_nt} \;,
\end{equation}
with nodes and weights, respectively,
\begin{equation}
x_n \equiv 1 + r + 2\sqrt{r}\cos(\pi\tau_n) \;, \quad
w_n \equiv \frac{\mu_1\sqrt{r}}{L}{\cdot}\frac{\sin^2(\pi\tau_n)}{\alpha + \cos(\pi\tau_n)}
     {\cdot}R_{\text{cut}}(-x_n) \;.
\end{equation}
Quadrature for the PDF is given by
\begin{equation}
P_{\text{cut}}(t) \simeq \sum_{n=1}^L w_nx_n e^{-x_nt} \;.
\end{equation}

Equation (\ref{FcutV}) may also be expressed in a form suitable for Gauss-Chebyshev
quadrature of the first kind, namely,
\begin{equation}
\bar{F}_{\text{cut}}(t) = \frac{\mu_1\sqrt{r}}{\pi}e^{-(r+1)t}\int_{-1}^{+1} \frac{dv}{\sqrt{1 - v^2}}\,
     e^{-2\sqrt{r}tv}{\cdot}\frac{(1-v)R_{\text{cut}}(-2\sqrt{r}(v+\alpha))}{1 + (\alpha - 1)/(1 + v)} \;.
\label{FcutVFirst}
\end{equation}
It is known that Gauss-Chebyshev quadrature of the first kind is equivalent
to the mid-point rule on the unit circle, which gives rise to the same quadrature scheme as
described for the trapezoidal rule, but where the nodes are generated by
\mbox{$\tau_n = (n - 1/2)/L$}
for
\mbox{$n = 1,2,\ldots,L$},
so that
\mbox{$0 < \tau_n < 1$}
which avoids any potential singularities.
We find that the most numerically robust scheme over the entire range of traffic intensities $r$
corresponds to the mid-point rule with the quadrature weights expressed in the form
\begin{equation}
w_n = \frac{2\mu_1\sqrt{r}}{L}{\cdot}R_{\text{cut}}(-x_n){\cdot}\frac{\sin^2(\pi\tau_n/2)}
     {1 + (\alpha-1)/[2\cos^2(\pi\tau_n/2)]} \;.
\end{equation}

The trapezoidal rule is generally a poorly performing quadrature scheme.
However, it has been shown to be exponentially convergent for some specific classes of integrand,
such as those comprising smooth periodic functions in which the integration
extends over a period \citep{BP:Rice73, BP:Trefethen14}.
It should be noted that the integrand above is unit-periodic in $\tau$.
Because of this, and the fact that it is equivalent to Gauss-Chebyshev quadrature, the
trapezoidal rule is an efficient way to compute the BP distribution. This is also true for the
mid-point rule.

Alternatively, in order to attain a guaranteed level of accuracy, we may apply
an iterative trapezoidal or mid-point rule to the integral
\begin{equation}
\bar{F}_{\text{cut}}(t) = \int_0^1 d\tau\, w(\tau) e^{-x(\tau)t} \;,
\end{equation}
with
\begin{equation}
w(\tau) \equiv \mu_1\sqrt{r}\frac{\sin^2(\pi\tau)}{\alpha + \cos(\pi\tau)}
     {\cdot}R_{\text{cut}}(-x(\tau)) \;,
\end{equation}
and
\begin{equation}
x(\tau) = 1 + r + 2\sqrt{r}\cos(\pi\tau) \;.
\end{equation}
Gaussian quadrature schemes do not generally allow nodes and weights from a previous iteration
to be re-used when sequentially moving to a finer grid in order to improve accuracy on the way
to satisfying some predetermined convergence criterion.
However, the trapezoidal rule can be iterated by doubling the number of grid points at each step.
Only half the grid points need to be computed, as the grid points from the previous iteration supply
the other half \citep{BP:Press07}.
Similar considerations apply to the iterative mid-point rule, but the grid size needs to be tripled at
each step, with a third of the total number of points coming from the previous iteration.
The large-$t$ asymptotics are derived from (\ref{FcutU}) by
\begin{align}
\begin{aligned}
\bar{F}_{\text{cut}}(t) &\asym{t\to\infty} \frac{2\mu_1\sqrt{r}}{\pi a}R_{\text{cut}}(-4a\sqrt{r})
     e^{-4a\sqrt{r}t}\int_0^\infty du\,\sqrt{u} e^{-4\sqrt{r}tu} \\
&= \frac{\mu_1}{8\sqrt{\pi}ar^{1/4}}R_{\text{cut}}(-4a\sqrt{r})t^{-3/2}e^{-4a\sqrt{r}t} \;.
\end{aligned}
\end{align}

We shall now look at the pole contribution to the BP distribution.
Poles occur in (\ref{PhiG}) when
\mbox{$\psi_+(s) - G_{22}(s) = 0$}.
Let us suppose that these are found at
\mbox{$s = \zeta_\ell$},
\mbox{$\ell = 1,2,\ldots$}.
Since the first term
\mbox{$G_{11}(s)$}
will not contribute to the residue,
the pole components of the MGF are given by the residues
\begin{equation}
\phi_1(s)/\mu_1 \underset{\text{pole}}{\longrightarrow} \frac{G_{12}^2(\zeta_\ell)}
     {\psi'_+(\zeta_\ell) - G'_{22}(\zeta_\ell)} \;,
\end{equation}
where the prime denotes differentiation.
Since, for general $s$,
\mbox{$\psi'_+(s) = \psi_+(s)/\sqrt{\Delta(s)}$},
we have at a pole
\mbox{$\psi'_+(\zeta_\ell) = G_{22}(\zeta_\ell)/\sqrt{\Delta(\zeta_\ell)}$}.
Then, on setting
\begin{equation}
J_{22}(s) \equiv \Pi^2(s)G'_{22}(s) = \sum_{k=1}^{N-1}v^{(k)}_{22}
     \prod_{\ell = 1 \atop \ell\neq k}^{N-1} (s - \xi_\ell)^2 \;,
\end{equation}
we obtain
\begin{equation}
\phi_1(s)/\mu_1 \underset{\text{pole}}{\longrightarrow} \frac{H^2_{12}(\zeta_\ell)\sqrt{\Delta(\zeta_\ell)}}
     {\Pi(\zeta_\ell)H_{22}(\zeta_\ell) + J_{22}(\zeta_\ell)\sqrt{\Delta(\zeta_\ell)}} \;.
\end{equation}
Therefore, the general result for the pole contribution to the BP SF is
\begin{equation}
\bar{F}_{\text{pol}}(t) = \mu_1\sum_{\ell = 1}^{n_{\text{pol}}}
     J(\zeta_\ell)e^{\zeta_\ell t}/(-\zeta_\ell) \;,
\end{equation}
with
\begin{equation}
J(s) \equiv \frac{H_{12}^2(s)\sqrt{\Delta(s)}}{\Pi(s)H_{22}(s) + J_{22}(s)\sqrt{\Delta(s)}} \;,
\label{JFn}
\end{equation}
and
\begin{equation}
\psi_+(\zeta_\ell) - G_{22}(\zeta_\ell) = 0 \;.
\end{equation}
The number of $\zeta$-poles can be determined {\it a priori}\ according to
\begin{equation}
n_{\text{pol}} = \sum_{\ell=1}^{N-1}I[-\xi_\ell < x_-] - I[G_{12}(0) \leq 0]
     + I[G_{22}(-x_-) - \psi_+(-x_-) \ge 0)] \;.
\end{equation}
The notation
\mbox{$I[\mathcal{A}]$}
is for the indicator function such that
\mbox{$I[\mathcal{A}] = 1$}
if $\mathcal{A}$ is true, and is zero otherwise.
The initial sum counts the contributing $\xi$-asymptotes, while the next two
terms account for boundary effects.
The maximum number of poles that can contribute is $N-1$, and this occurs when
\mbox{$r \to 0^+$},
in which case $x_\pm = 1$, so that the cut disappears.
Figure~\ref{BPPoles} displays the generation of the $\zeta$-poles as the intersection of a
sideways half-parabola defined by $\psi_+(s)$ with a tan-like function generated by the
$\xi$-poles of $G_{22}(s)$.
The tip of the parabola lies on the near edge of the cut.
Thus, the parabola can only intercept the $\xi$-pole curves
at points between the origin and the near cut boundary.
The full BP distribution is given by the sum of the pole and cut contributions, {\it i.e.}\
\mbox{$\bar{F}(t) = \bar{F}_{\text{pol}}(t) + \bar{F}_{\text{cut}}(t)$}
in the case of the SF.

By differentiating $\phi(s)$ at
\mbox{$s = 0$},
we can obtain an expression for the mean BP given by
\begin{equation}
\langle T_{\text{bp}}\rangle_{T_{\text{bp}}} = -\mu_1\left[G'_{11}(0) +
     \frac{2G_{12}(0)G'_{12}(0)}{1/r - G_{22}(0)}
     - G^2_{12}(0)\frac{1/(r(1-r)) - G'_{22}(0)}{\left(1/r - G_{22}(0)\right)^2}
     \right] \;.
\end{equation}
This provides a useful numerical diagnostic when compared with the exact result \citep{BP:Artalejo01}
\begin{equation}
\langle T_{\text{bp}}\rangle_{T_{\text{bp}}} = \frac{1}{r}\left[\frac{(Nr)^N}{N!(1-r)} +
     \sum_{k=1}^{N-1} \frac{(Nr)^k}{k!}\right]
     = \frac{1}{r}\left[\frac{(Nr)^N}{N!(1-r)} + e^{Nr}{\cdot}\frac{\Gamma(N,Nr)}{\Gamma(N)} - 1\right] \;.
\label{BPMeanEx}
\end{equation}
Here,
\mbox{$\Gamma(n, x)$}
denotes the upper incomplete gamma function with shape parameter $n$.

The upshot of the foregoing analysis is that the BP distribution is arbitrarily well represented as
an exponential mixture.
Thus, suppose that the RV $X$ is an exponential mixture, so that
\begin{equation}
\bar{F}_X(x) = \sum_{\ell=1}^L w_\ell e^{-x/u_\ell}
\end{equation}
for the SF.
Then, we have for the mean
\begin{equation}
\langle X\rangle_X \equiv \int_0^\infty dx\, xP(x)
= \sum_{\ell=1}^L u_\ell w_\ell \;.
\end{equation}
For the log-mean,
\begin{equation}
\langle\ln X\rangle_X \equiv \int_0^\infty dx\, P(x)\ln x
= \sum_{\ell=1}^L w_\ell\ln u_\ell - \gamma_{\text{e}} \;,
\end{equation}
where
\mbox{$\gamma_{\text{e}} \equiv -\psi(1) \simeq 0.5772$}
denotes Euler's constant.
For the differential entropy,
\begin{align}
\begin{aligned}
\langle -\ln P(X)\rangle_X &\equiv -\int_0^\infty dx\, P(x)\ln P(x) \\
&= -\int_0^\infty dx\, e^{-x} \sum_{\ell=1}^L w_\ell \ln\left[
     \sum_{k=1}^L \left(w_k/u_k\right) e^{-x(u_\ell/u_k)}\right] \;,
\end{aligned}
\end{align}
which is easily computed via Gauss-Laguerre quadrature.
In the present application to the BP, the nodes should account for the
discrete poles as well as the cut.
Thus,
\mbox{$L = L_{\text{cut}} + L_{\text{pol}}$}.
The discrete poles contribute nodes
\mbox{$x_\ell = -\zeta_\ell$}
with weights
\mbox{$w_\ell = -\mu_1 J(\zeta_\ell)/\zeta_\ell$}.
From the requirement that
\mbox{$\bar{F}_X(0) = 1$},
we have that
\mbox{$\sum_{\ell = 1}^Lw_\ell = 1$}.
Thus, the node/weight pairs
\mbox{$(u_\ell,w_\ell)$}
may be thought of as describing the empirical SF for some RV
\mbox{$0 \leq U < \infty$}
according to
\begin{equation}
\bar{F}_{U}(u) = \sum_{\ell = 1}^L w_\ell I(u_\ell > u) \,,
\label{EmpTxtSF}
\end{equation}
which we may refer to as the texture distribution by analogy with compound Gaussian
clutter distributions in radar detection theory \citep{BP:Rosen22}.
This compound representation implies that the BP RV has the product form
\mbox{$T_{\text{bp}} = UV$},
where here $V$ is a unit-mean exponentially distributed RV,
which provides a very convenient scheme for generating random numbers drawn from the
BP distribution.

\begin{figure}
\FIGURE
{\includegraphics[width=\wscl\linewidth, height=\hscl\linewidth]{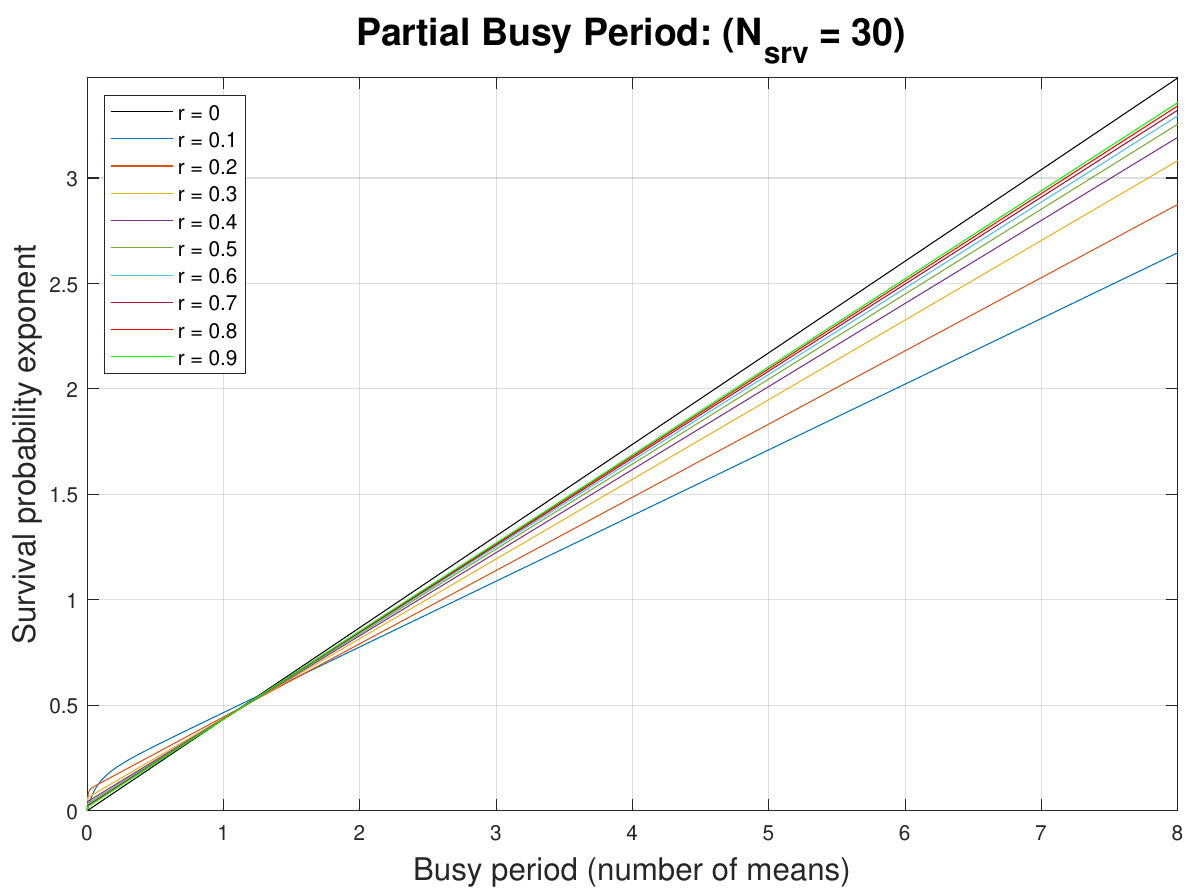}}
{\hphantom{x}\label{BPSFTail}}
{BP survival function for various traffic intensities ($r$)
     and $N_{\text{srv}} = 30$ servers, showing long-time behaviour.
     The negative base-10 logarithm of the survival function is plotted.}
\end{figure}

\begin{figure}
\FIGURE
{\includegraphics[width=\wscl\linewidth, height=\hscl\linewidth]{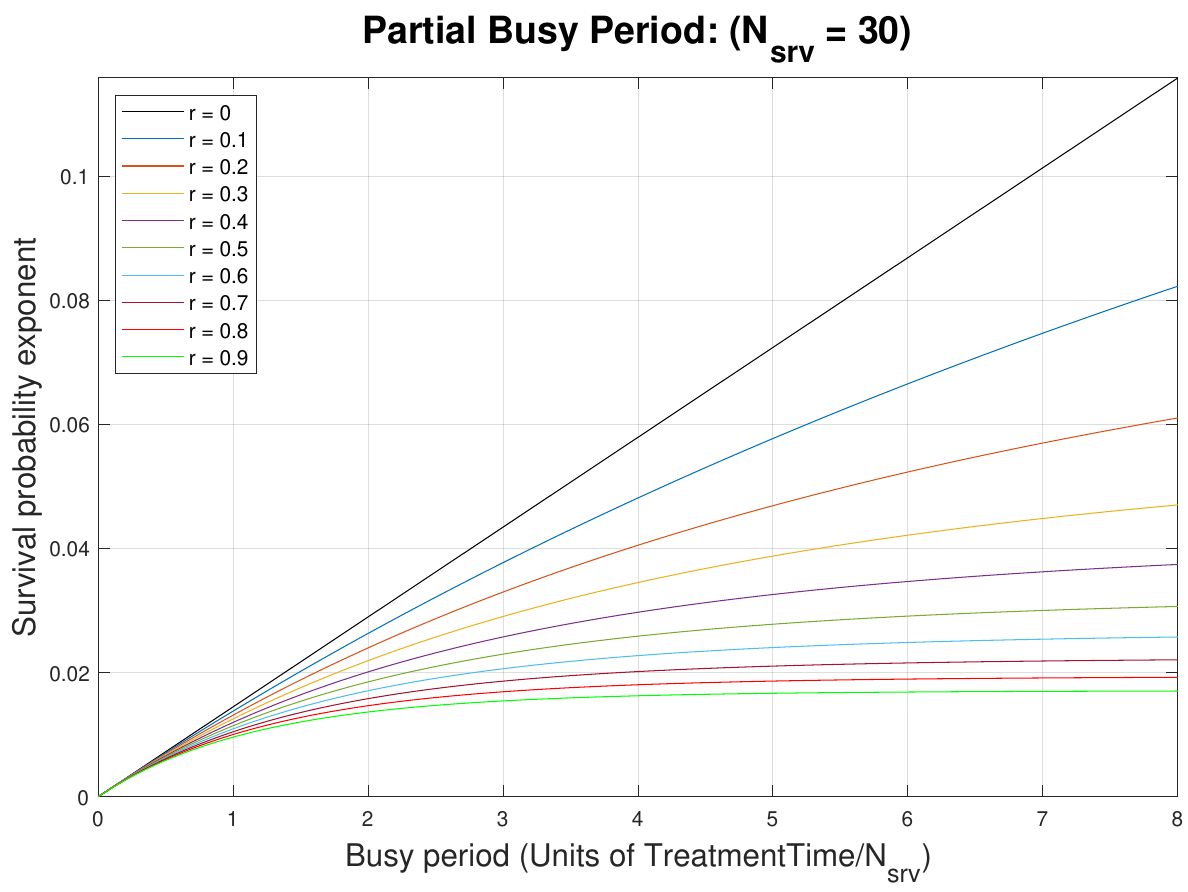}}
{\hphantom{x}\label{BPSFTailZoom}}
{BP survival function for various traffic intensities ($r$)
     and $N_{\text{srv}} = 30$ servers, showing short-time behaviour.
     The negative base-10 logarithm of the survival function is plotted.}
\end{figure}

Figures~\ref{BPSFTail} and \ref{BPSFTailZoom} display the logarithm of the BP SF for the case of
\mbox{$N = 30$}
servers and a variety of traffic intensities.
In Figure \ref{BPSFTailZoom}, the horizontal axis records the (dimensionless) time as directly given by
\mbox{$\bar{F}(t) = \bar{F}_{\text{pol}}(t) + \bar{F}_{\text{cut}}(t)$},
in which case one is observing the short-time behaviour of the distribution.
Note that the MGF for the problem formulated in terms of physical units
can be recovered from the dimensionless problem via
\mbox{$\phi_{\text{phys}}(s) = \phi_1(s/s_{\text{scl}})$}
with frequency scale
\mbox{$s_{\text{scl}} = N\mu$}.
For the SF, this implies
\mbox{$\bar{F}(t) = \bar{F}_{\text{phys}}(t_{\text{scl}}{\cdot}t)$}
with timescale
\mbox{$t_{\text{scl}} = 1/s_{\text{scl}} = 1/(N\mu)$}.
In Figure \ref{BPSFTail}, the horizontal axis measures time in units of the mean BP.
Thus, the function
\mbox{$\bar{F}(m_{\text{bp}}t)$}
is being plotted, and this illustrates the quite distinct long-time behaviour of the distribution.

Since the regeneration cycle time is the sum of a BP and an idle period, and these are mutually
independent, it follows that the regeneration cycle distribution is the convolution of the BP
distribution and the inter-arrival distribution. Therefore, the SF for the regeneration cycle time can
formally be expressed as
\begin{equation}
\bar{F}_{\text{reg}}(t) \simeq \sum_{\ell=1}^L w_\ell{\cdot}
     \frac{r e^{-x_\ell t} - x_\ell e^{-rt}}{r - x_\ell} \;.
\end{equation}
This expression should be numerically unproblematic for
\mbox{$r < 1/4$}. Above this value, the inter-arrival pole at
\mbox{$s = -r$}
lies within the cut
\mbox{$\left\{s: (1-\sqrt{r})^2 \leq -s \leq (1+\sqrt{r})^2\right\}$}.
It should be noted, however, that for appreciable values of the traffic intensity,
the idle period gives a negligible contribution to the regeneration cycle.

\subsection{Boundary Case: $r \to 0$}
Numerical instabilities arise for small values of the total traffic intensity $r$.
However, in the
\mbox{$r\to 0^+$}
limit, the $v_{12}$ become negligible, in which case we have
\begin{equation}
\phi_1(s)/\mu_1 \simeq \sum_{k=1}^{N-1}\frac{|v^{(k)}_1|^2}{s - \xi_k}
     = \sum_{k=1}^{N-1}\frac{v_{11}^{(k)}}{s - \xi_k} \;.
\end{equation}
It follows, for the BP PDF and SD, respectively, that
\begin{equation}
P(t) \simeq \mu_1\sum_{k=1}^{N-1} v_{11}^{(k)}e^{\xi_k t} \;, \quad
\bar{F}(t) \simeq \mu_1\sum_{k=1}^{N-1} \frac{v_{11}^{(k)}}{|\xi_k|}e^{\xi_k t} \;,
\end{equation}
recalling that
\mbox{$\xi_k < 0$}.
A criterion can be developed to signal when to use this approximation.
We have found that the test
\mbox{$\max\bigl\{|v_{12}^{(k)}|: k = 1,2,\ldots,N-1\bigr\} < 10^{-4}$}
works adequately.

\subsection{Boundary Case: $r \to 1^-$}
Poles
\mbox{$s = \zeta_\ell$}
occur when
\mbox{$\psi_+(s) - G_{22}(s) = 0$}.
When the smallest pole is very close to zero, numerical solution of this equation
can present difficulties.
In this case, we can linearize about
\mbox{$s = 0$}
to obtain
\begin{align}
\begin{aligned}
\psi_+(s) &= \psi_+(0) + s\psi'_+(0) + O(s^2) \\
     &= \frac{1}{r} + \frac{1}{r(1-r)}s + O(s^2) \;, \\
G_{22}(s) &= G_{22}(0) + sG'_{22}(0) + O(s^2) \;,
\end{aligned}
\end{align}
which leads to the approximate solution
\begin{equation}
s \simeq -\frac{1 - rG_{22}(0)}{1/(1-r) - rG'_{22}(0)}
     = -\frac{\mu_1}{\sqrt{\Lambda}}{\cdot}\frac{G_{12}(0)}{1/(1-r) - rG'_{22}(0)} \;.
\label{PoleNearZero}
\end{equation}
In order to catch bad numerical behaviour due to arithmetic underflow,
when the condition
\mbox{$\left|1/r - G_{22}(0)\right| < 10^{-9}$}
is detected, the sign check
\begin{equation}
\sgn(G_{12}(0)) = \sgn(1/ r- G_{22}(0))
\end{equation}
is performed.
If this check succeeds, then $\zeta_1$ is replaced by the value obtained from the latter of (\ref{PoleNearZero}).
On the other hand, failure of the test indicates that the lowest pole has missed altogether, in which case
$\zeta_1$ obtained from (\ref{PoleNearZero}) is appended to the collection of poles $\zeta_k$.

A variety of diagnostic tests is available to track the validity of the numerical computations.
For example, it can be checked (i) whether the expected number of $\zeta$-poles is being generated,
(ii) whether equation (\ref{BPMeanEx}) for the mean holds to some desired accuracy,
(iii) whether the relationship (\ref{G12G22}) holds to some desired accuracy and, among others,
(iv) whether the MGF as given by (\ref{PhiG}) evaluates to unity at the origin.

\subsection{Boundary Case: $r = 1$}
At the boundary of the ergodic region, where
\mbox{$r = 1$},
there is no pole contribution.
Hence, the BP SF comes entirely from the cut, and reads
\begin{align}
\begin{aligned}
\bar{F}(t) &= \frac{2\mu_1}{\pi}\int_0^1 \frac{du}{u}\, e^{-4tu}R_{\text{cut}}(-4u)
     {\cdot}\sqrt{u(1-u)} \\
&\asym{t\to\infty} \frac{\mu_1}{\pi\sqrt{t}}R_{\text{cut}}(0)\int_0^{4t} \frac{du}{\sqrt{u}}\, e^{-u} \\
&\asym{t\to\infty} \frac{\mu_1}{\pi\sqrt{t}}R_{\text{cut}}(0)\int_0^\infty \frac{du}{\sqrt{u}}\, e^{-u} \\
&= \mu_1 R_{\text{cut}}(0){\cdot}\frac{1}{\sqrt{\pi t}} \;,
\end{aligned}
\end{align}
where we may recall that, for
\mbox{$r = 1$},
\mbox{$\mu_1 R_{\text{cut}}(0) = 1/\prod_{n=1}^{N-1}\mu_n$}.
We see that the BP distribution becomes power-law tailed in the non-ergodic region.
We can write
\begin{equation}
\bar{F}(t) \asym{t\to\infty} \frac{1}{1 + \sqrt{\pi t}/[\mu_1 R_{\text{cut}}(0)]}
\end{equation}
to obtain a form that is unity at the origin, rather than diverging.
One can also adopt the Lomax form
\begin{equation}
\bar{F}(t) \asym{t\to\infty} \frac{1}{\left(1 + \tfrac{\pi}{[\mu_1 R_{\text{cut}}(0)]^2}t\right)^{1/2}} \;.
\label{SFLomax}
\end{equation}
For
\mbox{$N = 1,2$},
the Lomax form gives a close fit across the whole domain.

\section{Algebraic Approach}
\label{Algebraic}
In this section, we exploit
the tridiagonal structure of the problem to derive a continued-fraction representation of the
BP MGF that can be used to yield explicit expressions on a case-by-case basis for small numbers of servers.
We then proceed to analyse the underlying tridiagonal matrix in a different way that allows explicit
representation of the cut function for any number of servers in terms of a family of polynomials
(to which we refer as the cut polynomials). Finally, we derive an easily solvable two-dimensional
system of recurrence relations that enables
determination of the coefficients of the cut polynomials which, in turn, allows for the calculation
of their roots. As demonstrated in a subsequent section, these roots provide complete information for the
construction of closed-form expressions for the BP distribution functions.

\subsection{Continued Fraction}
Let
\mbox{$M_0 \equiv \diag[s,s,\cdots,s-r\psi_-(s)] - M$},
where $M$ is the full tridiagonal birth-death process generator matrix for the BP problem
as given in (\ref{TriDiagM}), and let $M_1$
be the matrix obtained from $M_0$ by omitting the first row and column.
We then define the general matrix $M_k$ recursively as the one obtained by
omitting the first row and column from $M_{k-1}$.
This can be continued until one reaches the scalar
\begin{equation}
M_{N-2} = s + \mu_{N-1} - r \psi_-(s) \;.
\end{equation}
Cramer's rule yields
\begin{equation}
\left(M_0^{-1}\right)_{11} = \det M_1/\det M_0 \;,
\end{equation}
so that
\begin{equation}
\phi_1(s) = \left(M_0^{-1}\right)_{11}\mu_1 = \mu_1\det M_1/\det M_0 \;.
\end{equation}
By applying the cofactor expansion in the first row,
we are led to the recurrence relation
\begin{equation}
\det M_k = \sigma_{k+1}\det M_{k+1} - r\mu_{k+2}\det M_{k+2} \;,
\end{equation}
subject to
\mbox{$M_{N-2} \equiv \sigma_{N-1} - r\psi_-$},
\mbox{$M_{N-1} \equiv 1$},
where we have introduced the notation
\mbox{$\sigma_k \equiv s + r + \mu_k$}.
By writing
\begin{equation}
\phi_1(s) = \frac{\mu_1}{\sigma_1 - r\mu_2\det M_2/\det M_1} \;,
\end{equation}
on applying the recurrence once, and then successively applying the recurrence, we arrive at the
continued fraction representation
\begin{equation}
\phi_1(s) = \frac{1}{r}\left[\frac{r\mu_1}{\sigma_1 - {\cdot}}, \frac{r\mu_2}{\sigma_2 - {\cdot}}, \cdots,
   \frac{r\mu_{N-1}}{\sigma_{N-1} - r\psi_-}\right] \;.
\label{cf}
\end{equation}
One may note that if we set
\mbox{$\eta_N(s) \equiv r\psi_-(s)$},
and
\begin{equation}
\eta_k(s) = \frac{r\mu_k}{\sigma_k - \eta_{k+1}(s)} \;,
\end{equation}
for
\mbox{$k = N-1, N-2, \ldots, 1$},
then we have
\mbox{$\phi_1(s) = \eta_1(s)/r$}.

Setting
\mbox{$b(s) \equiv s + r + 1$},
we can write
\begin{equation}
\psi_\pm(s) = \frac{1}{r}\left[b(s) \pm \sqrt{b^2(s) - 4r}\right] \;.
\end{equation}
Since
\begin{equation}
r\psi_\pm^2(s) - b(s)\psi_\pm(s) + 1 = 0 \;,
\end{equation}
it follows that
\begin{equation}
r\psi_+(s) + r\psi_-(s) = b(s) \;, \quad
r\psi_+(s) \cdot r\psi_-(s) =  r \;,
\end{equation}
from which we also obtain the identity
\begin{equation}
(\sigma_k - r\psi_-)(\sigma_k - r\psi_+) = (\mu_k - 1)\sigma_k + r \;.
\end{equation}
With the aid of these relationships,
(\ref{cf}) leads directly to the following explicit forms linear in
\mbox{$\psi_-(s)$}:
For
\mbox{$N = 2$},
with
\mbox{$\mu_k = k/2$},
\begin{equation}
\phi_1(s) = \mu_1\frac{\mu_1 - 1 + r\psi_-(s)}{(\mu_1-1)\sigma_1 + r} \;.
\end{equation}
For
\mbox{$N = 3$},
with
\mbox{$\mu_k = k/3$},
\begin{equation}
\phi_1(s) = \frac{\mu_1\sigma_1^2 - (r - 1 + 2\mu_1^2)\sigma_1
     + 2\mu_1(\mu_1 r - \mu_1 + 1) - 2\mu_1r\psi_-(s)}
     {\sigma_1^3 + (\mu_1 - r)\sigma_1^2 - 4r^2\mu_1} \;.
\end{equation}
For
\mbox{$N = 1$},
we trivially have
\mbox{$\phi_1(s) = \psi_-(s)$}.
While the continued-fraction representation can be used to derived explicit expressions
for the MGF when the number of servers is small, there does not seem to be any obvious way
to extract an explicit closed-form result for the general case.
A finite continued fraction representation was previously considered by \citet{BP:Daley98},
who make the same observation that `there is no simple non-recursive closed-form solution'.
Prior to this, \citet{BP:Conolly74} had also derived a continued fraction representation.

\subsection{Cut Polynomials}
Further progress for the general case can be made by observing the general structure
\begin{equation}
\phi_1(s) = \mu_1\frac{N_0(\sigma_{N-1} - r\psi_-) + N_1}{D_0(\sigma_{N-1} - r\psi_-) + D_1} \;,
\end{equation}
where
\mbox{$N_{0,1}$}
and
\mbox{$D_{0,1}$}
are polynomials in $s$.
This representation follows from applying the cofactor expansion for
\mbox{$\det M_0$}
and
\mbox{$\det M_1$}
to the last rather than the first row.
At a pole, which occurs when
\mbox{$\sigma_{N-1} - r\psi_- = -D_1/D_0$},
the numerator becomes rational:
\begin{equation}
N_0(\sigma_{N-1} - r\psi_-) + N_1 = \frac{N_1D_0 - N_0D_1}{D_0} \;.
\end{equation}
It is convenient to introduce the functions
\begin{align}
\begin{aligned}
N_\pm &\equiv N_0(\sigma_{N-1} - r\psi_\pm) + N_1 \;, \\
D_\pm &\equiv D_0(\sigma_{N-1} - r\psi_\pm) + D_1 \;, \\
\end{aligned}
\end{align}
so that
\mbox{$\phi_1(s) = \mu_1 N_-(s)/D_-(s)$}.
One can show that
\begin{equation}
D_+D_- = D_1^2 + D_0\left[((\mu_{N-1}-1)\sigma_{N-1} + r)D_0
     + (\sigma_{N-1} + \mu_{N-1} - 1)D_1\right] \;,
\end{equation}
which is a polynomial in $s$.
We also have
\begin{equation}
N_-D_+ = A + (N_1D_0 - N_0D_1){\cdot}r\psi_- \;,
\end{equation}
for some polynomial $A(s)$,
and one can show that
\begin{equation}
N_1D_0 - N_0D_1 = \frac{1}{r\mu_1}\prod_{k=1}^{N-1}(r\mu_k) \;.
\end{equation}
Thus, we see that the representation
\begin{equation}
\phi_1(s) = \mu_1 \frac{N_-(s)D_+(s)}{D_+(s)D_-(s)}
\end{equation}
has polynomial denominator and is linear in $\psi_-(s)$.

On the cut in the $s$-plane,
that is due to square-root term in
\mbox{$\psi_-(s) = [b(s) - \sqrt{|\Delta(s)|}]/(2r)$},
the function $\phi_1(s)$
contributes only the component
\begin{equation}
\phi_{\text{cut}}(s) = (\pm i)\frac{\prod_{k=1}^{N-1}(r\mu_k)}{D_+(s)D_-(S)}{\cdot}
     \frac{\sqrt{|\Delta(s)|}}{r} \;,
\end{equation}
where
\mbox{$|\Delta(s)| = 4r - b^2(s)$}
on the cut.
This allows us to identify the cut function introduced previously as
\begin{equation}
R_{\text{cut}}(s) = \frac{1}{r\mu_1}{\cdot}\frac{\prod_{k=1}^{N-1}(r\mu_k)}{D_+(s)D_-(s)}
     = -\frac{1}{r\mu_1^2}{\cdot}\frac{\prod_{k=1}^{N-1}(r\mu_k)}{C_N\left(\sigma(s)\right)} \;,
\label{FcutC}
\end{equation}
where we find it convenient to introduce the cut polynomials $C(\sigma)$,
defined as functions of the variable
\mbox{$\sigma \equiv \sigma_1 = s + r + \mu_1$}
via
\mbox{$C_N\left(\sigma(s)\right) \equiv -\mu_1 D_+(s)D_-(s)$}.
The first few of these are given as follows:
For
\mbox{$N = 2$},
\mbox{$C_2(\sigma) = \sigma - r/\mu_1$}.
For
\mbox{$N = 3$},
\begin{equation}
C_3(\sigma) = \sigma^3 + (\mu_1 - r)\sigma^2 - 4r^2\mu_1 \;.
\end{equation}
For
\mbox{$N = 4$},
\begin{equation}
\begin{split}
C_4(\sigma)& = \sigma^5 + \left(1 - r\right)\sigma^4  + \left(5\mu_1 - 6r\right)\mu_1\sigma^3 \\
&    {}+ \left(2\mu_1^2 - 13\mu_1r + r^2\right)\mu_1\sigma^2
     - 2\left(1 - 7r\right)\mu_1^2r\sigma
     + 2\mu_1^2r^2 - r^3 \;.
\end{split}
\label{C4}
\end{equation}
In general,
\mbox{$C_N(\sigma) = \sigma^{2N-3} + \cdots$},
provided
\mbox{$N > 1$}.
For consistency with (\ref{FcutC}), one may set
\mbox{$C_1(\sigma) \equiv -1/\mu_1^2$}
for the single server case.

The residue of $\phi_1(s)$ at a given pole $s$,
defined by
\begin{equation}
\phi_{\text{res}}(s) \equiv \lim_{s' \to s}(s' - s)\phi_1(s') \;,
\end{equation}
is
\begin{equation}
\phi_{\text{res}}(s) = \mu_1\frac{N_1D_0 - N_0D_1}{D_0}{\cdot}
     \frac{D_+}{(D_+D_-)'} \;,
\end{equation}
where the prime in
\mbox{$(D_+D_-)'$}
denotes that the pole factor has been divided out of the polynomial, {\it i.e.}\
\begin{equation}
(D_+D_-)'(s) \equiv \lim_{s' \to s} N_-(s')D_+(s')/(s' - s) \;.
\end{equation}
We can derive the general relationship
\begin{equation}
D_+ + D_- = D_0(\sigma_{N-1} + \mu_{N-1} - 1) + 2D_1 \;,
\end{equation}
and note that
\mbox{$D_- = 0$}
at a pole, to obtain that, at a pole,
\begin{equation}
D_+ = D_0(\sigma_{N-1} + \mu_{N-1} - 1) + 2D_1 \;.
\end{equation}
Consequently, at a pole $s$, the residue term is
\begin{equation}
\phi_{\text{res}}(s) = \frac{1}{r}\prod_{k=1}^{N-1}(r\mu_k){\cdot}
     \frac{\sigma_{N-1} + \mu_{N-1} - 1 + 2D_1/D_0}{(D_+D_-)'} \;,
\end{equation}
It may be observed that the polynomials $N_0,N_1$ have completely dropped out of the calculations.
The full BP distribution is obtained from the sum of the pole and cut contributions to the MGF.
For the PDF,
\begin{equation}
P(t) = \sum_{\ell\in\text{poles}}\phi_{\text{res}}(\zeta_\ell)e^{\zeta_\ell t}
     + \int_{\text{cut}}\frac{ds}{2\pi}\, \phi_{\text{cut}}(s) e^{st} \;.
\end{equation}

\subsection{Recurrence Relations}
Let $M^{(0)}_N$ denote the the tridiagonal matrix with main diagonal
\mbox{$(\sigma_1, \sigma_2, \cdots, \sigma_{N-1})$},
lower sub-diagonal
\mbox{$(-\mu_2, -\mu_3, \cdots, -\mu_{N-1})$}
and constant upper sub-diagonal with common value $-r$.
Then,
\begin{equation}
D_0 = \det M^{(0)}_{N-2} \;, \quad
D_1 = -r\mu_{N-1}\det M^{(0)}_{N-3} \;.
\end{equation}
Let $M^{(1)}_N$ denote the the tridiagonal matrix with main diagonal
\mbox{$(\sigma_2, \sigma_3, \cdots, \sigma_{N-1})$},
lower sub-diagonal
\mbox{$(-\mu_3, -\mu_4, \cdots, -\mu_{N-1})$}
and constant upper sub-diagonal with common value $-r$.
Then,
\begin{equation}
N_0 = \det M^{(1)}_{N-2} \;, \quad
N_1 = -r\mu_{N-1}\det M^{(1)}_{N-3} \;.
\end{equation}
In terms of the matrix $M$ in (\ref{TriDiagM}),
\mbox{$M^{(0)}_N = s\mathbb{I} - M$}
and $M^{(1)}_N$ is $M^{(0)}_N$ with the first row and column removed.
We have the recurrences
\begin{equation}
\det M^{(j)}_{N} = \sigma_{N-1}\det M^{(j)}_{N-1} - r\mu_{N-1}\det M^{(j)}_{N-2} \;,
\end{equation}
for each
\mbox{$j = 0,1$}.
Only the
\mbox{$X_N \equiv \det M^{(0)}_{N}$}
are significant for general computation of the MGF.
Setting
\mbox{$\sigma_n = \sigma_1 + \hat{\mu}_n$}
with
\mbox{$\hat{\mu}_n \equiv \mu_n - \mu_1$},
we see that $X_N$ is obtained from the recurrence
\begin{equation}
X_n = (\sigma_1 + \hat{\mu}_{n-1})X_{n-1} - r\mu_{n-1}X_{n-2} \;,
\end{equation}
for
\mbox{$n = 2,3,\ldots,N$},
subject to
\mbox{$X_0 = 0$},
\mbox{$X_1 = 1$}.
To obtain an explicit expansion in powers of $\sigma_1$, we write
\begin{equation}
X_n = \sum_{k=0}^{n-1} x_n^{(k)}\sigma_1^k \;, \quad x_n^{(n-1)} \equiv 1 - \delta_{n0} \;, \quad x_n^{(-1)} \equiv 0 \;
\end{equation}
for
\mbox{$n = 0,1,2,\ldots$}.
This yields the simple two-dimensional recursion
\begin{equation}
x_n^{(\ell)} = \hat{\mu}_{n-1}x_{n-1}^{(\ell)} + x_{n-1}^{(\ell-1)} - r\mu_{n-1}x_{n-2}^{(\ell)} \;,
\end{equation}
for
\mbox{$0 \leq \ell \leq n-3$}.
It is straightforward to show that
\begin{equation}
x_n^{(n-2)} = \sum_{k=0}^{n-1} \hat{\mu}_k \;,
\end{equation}
for
\mbox{$n = 1,2,\ldots$},
which is required to seed the recursion.
Alternatively, one may achieve this by introducing
\mbox{$x_n^{(n)} \equiv 0$},
\mbox{$n \geq 0$}.
This scheme enables us, in principle, to calculate all the coefficients
of the general cut polynomial by invoking the identities
\mbox{$D_0 = X_{N-2}$},
\mbox{$D_1 = -r\mu_{N-1}X_{N-3}$}.
Table~\ref{tab:ND} provides expressions for $N_{0,1}$, $D_{0,1}$
for server numbers
\mbox{$N = 2,3,4$}.

\begin{table}
\renewcommand{\arraystretch}{1.3}
\TABLE
{MGF polynomials\label{tab:ND}}
{\begin{tabular}{|c||c|c|c|c|}
\hline
$N$ & $N_0$    & $N_1$       &  $D_0$                      & $D_1$ \\
\hline\hline
$2$ & $0$       & $1$        & $1$                         & $0$ \\
\hline
$3$ & $1$       & $0$        & $\sigma_1$                  & $-r\mu_2$ \\
\hline
$4$ & $\sigma_2$ & $-r\mu_3$ & $\sigma_1\sigma_2 - r\mu_2$ & $-r\mu_3\sigma_1$ \\
\hline
\end{tabular}}
{}
\end{table}

The critical value of the total traffic intensity is the value
\mbox{$r = r_{\text{c}}$}
at which the pole contribution disappears.
This occurs when
\mbox{$D_0(\sigma_{N-1} - \sqrt{r}) + D_1 = 0$}
at the cut extremity
\mbox{$s = -(1 - \sqrt{r})^2$}.
For
\mbox{$N = 2$},
it is
\mbox{$r_{\text{c}} = 1/2$}.
For
\mbox{$N = 3$},
\mbox{$y \equiv \sqrt{r_{\text{c}}}$}
solves the quadratic
\begin{equation}
2y^2 - 2y + \mu_1 = 0 \;,
\end{equation}
and for
\mbox{$N = 4$},
it solves the cubic
\begin{equation}
2y^3 - 45\mu_1^2 y^2 + y - 6\mu_1^3 = 0 \;.
\end{equation}
In Table~\ref{tab:rcrit},
we present numerical values for $r_{\text{c}}$
corresponding to the first few values of $N$.

\begin{table}
\renewcommand{\arraystretch}{1.3}
\TABLE
{Critical traffic intensities\label{tab:rcrit}}
{\begin{tabular}{|c||c|c|c|c|c|c|}
\hline
$N$ & $1$ & $2$ & $3$ & $4$ & $5$ & $6$ \\
\hline
$r_{\text{c}}$ & $0.0000$ & $0.5000$ & $0.6220$ & $0.8400$ & $0.9352$ & $0.9740$ \\
\hline
\end{tabular}}
{}
\end{table}

\section{Asymptotic Limits}
\label{Asym}
There are two natural ways to scale the problem in the limit of a large number of servers
\mbox{$N\to\infty$}.
The first is to keep the total traffic intensity $r$ constant as
\mbox{$N\to\infty$}.
The second is to keep the mean arrival and service rates
(namely $\lambda$, $\mu$, respectively) constant.
As either one of the two time-dimensional parameters simply serves to calibrate the
clock in the model, we can choose to set
\mbox{$\mu = 1$}.
This is equivalent to measuring model time in units of the mean service time,
and physical time can be recovered by appropriately rescaling the time axis of the resulting distributions.
Thus, the model essentially depend on only two parameters, which can be taken to be the number
of servers, $N$, and either $r$ or the partial traffic intensity
\mbox{$\rho \equiv \lambda/\mu$} \citep{BP:Gnedenko89}.
An advantage of the former parametrization is that the ergodic region of interest here
is characterized by a common bounded interval ({\it viz.}\
\mbox{$0 \leq r < 1$})
independent of $N$.

The limit of an infinite number of servers with given constant $\rho = r/N$
is referred to as the M/M/$\infty$ queue, and corresponds to our second scaling option.
It has been previously extensively studied by \citet{BP:Guillemin95},
and also earlier by \citet{BP:Morrison87}.
However, while this model is useful in some applications, it is also rather trivial
as a queueing model given that it
does not actually contain a queue.
We shall proceed to study the BP distributions of both scaling options.

Having explicit results for the
\mbox{$N\to\infty$}
limits, and comparing these with full calculations for finite given values of $N$,
enables one to decide, in creating and optimizing numerical algorithms, how far
one needs to go before simply being able to invoke asymptotic results.

\subsection{Constant-$\mathbf{\rho}$ Model}
To study this limit, it is convenient to consider the transpose problem.
Let
\mbox{$A \equiv s\mathbb{I} - M$},
where $M$ is given by (\ref{TriDiagM})\footnote{The normalization chosen in (\ref{PsiNorm}) is not appropriate here.
We take $\mu = 1$, rather than $\mu = 1/N$, so that $\lambda = \rho$ and adjust the matrix $M$ accordingly.}.
In the large-$N$ limit, $M$ becomes countably infinite dimensional
and the boundary condition for (\ref{PhiRecur}), {\it viz.}\
\mbox{$\phi_N = \psi_-(s)\phi_{N-1}$},
disappears as it gets pushed to infinity.
In Section~\ref{Spectral}, we considered a linear system which,
in the large-$N$ limit, now reads
\mbox{$A|\phi\rangle = |e_1\rangle$}
since we can set
\mbox{$|w\rangle = |e_1\rangle$},
noting that
\mbox{$\mu_1 = 1$}
here.
This can be inverted to yield the BP MGF
\mbox{$\phi_1(s) = \langle e_1|A^{-1}|e_1\rangle$}.
However, we have
\begin{equation}
\langle e_1|A^{-1}|e_1\rangle = \langle e_1|(A^{\sf T})^{-1}|e_1\rangle \;.
\end{equation}
Thus, we can also consider the transpose problem and compute
\mbox{$\varphi_1(s) \equiv \langle e_1|(A^{\sf T})^{-1}|e_1\rangle = \phi_1(s)$}.
Clearly
\mbox{$A^{\sf T} = s\mathbb{I} - M^{\sf T}$}
and, in the
\mbox{$N\to\infty$}
limit, the infinite tridiagonal matrix $M^{\sf T}$ is given by
\begin{equation}
M^{\sf T} =
\begin{bmatrix}
  -(\rho + \mu_1) & \mu_2           & 0               & \cdots     & \cdots          \\
  \rho            & -(\rho + \mu_2) & \mu_3           & \cdots     & \cdots          \\
  0               & \rho            & -(\rho + \mu_3) & \mu_4      & \cdots         \\
  \vdots          &                 & \ddots          & \ddots     & \ddots          \\
  \vdots          &                 &                 & \hphantom{-(\rho+\mu-3)}     & \ddots
\end{bmatrix} \;.
\label{TriDiagMInf}
\end{equation}
with
\mbox{$\mu_n \equiv n$}
being adopted here.
On setting
\mbox{$\varphi_0 \equiv 1/\rho$},
the implied recurrence relation for $\varphi_n$ reads
\begin{equation}
-\rho\varphi_{n-1} + (s + \rho + n)\varphi_n - (n+1)\varphi_{n+1} = 0 \;,
\label{VarPhiRecur}
\end{equation}
for
\mbox{$n = 1,2,\ldots$}.
Following \citep{BP:Guillemin95},
we next introduce the generating function
\begin{equation}
g(z) \equiv\sum_{n=1}^\infty z^{n-1}\varphi_n \;,
\end{equation}
so that
\mbox{$g(0) = \varphi_1 = \phi_1$}.
Summation over the recurrence (\ref{VarPhiRecur}) followed by some standard manipulations
gives rise to the differential equation
\begin{equation}
g'(z) + \left[\frac{s}{z-1} + \frac{1}{z} - \rho\right]g(z) = \frac{z - g(0)}{z(z-1)} \;.
\end{equation}
Solution via the integrating factor method yields
\begin{equation}
g(z) = \frac{e^{\rho z}}{(1 - z)^s}\int_0^1 d\xi\, e^{-\rho z\xi}(1 - z\xi)^{s-1}
     \left(g(0) - z\xi\right) \;.
\end{equation}
The vanishing of the residue at
\mbox{$z = 1$}
implies that
\begin{equation}
g(0){\cdot}\int_0^1 d\xi\, e^{-\rho\xi}(1 - \xi)^{s-1} =
     \int_0^1 d\xi\, e^{-\rho\xi}\xi(1 - \xi)^{s-1} \;.
\end{equation}
\citet{BP:Guillemin95} have defined the integrals
\begin{equation}
\mathcal{I}_\alpha(s,\rho) \equiv \int_0^1 d\xi\, e^{-\rho\xi}\xi^\alpha(1 - \xi)^{s} \;.
\end{equation}
In terms of these, we can write the result for the BP MGF as
\begin{equation}
\phi_1(s) = \mathcal{I}_1(s-1,\rho)/\mathcal{I}_0(s-1,\rho) \;.
\end{equation}
This result coincides with the MGF of \citep{BP:Guillemin95} for their congestion duration (denoted $\theta$)
in the case $C = 0$. These authors, however, did not proceed to recover the distribution from the MGF.
We shall also introduce the functions
\begin{equation}
\phi_1^{(\alpha)}(s) \equiv \mathcal{I}_{\alpha+1}(s-1,\rho)/\mathcal{I}_\alpha(s-1,\rho) \;,
\end{equation}
so that
\mbox{$\phi_1(s) = \phi_1^{(0)}(s)$}.
We note from the defining matrix equation
\mbox{$A|\phi\rangle = |e_1\rangle$},
that $\phi_1(s)$ has poles at the eigenvalues of the matrix $M^{\sf T}$.
Similarly, $\phi_1^{(1)}(s)$ has poles at the eigenvalues of the truncated problem
obtained by removing the first row and column from the matrix $M^{\sf T}$,
and these correspond to the zeros of $\mathcal{I}_1(s-1,\rho)$.

At this point, it is useful to recall that the Kummer hypergeometric function \citep{BP:Pearson17},
also known as the hypergeometric function of the first kind,
\mbox{$M(a,b,z) = \tensor[_1]{F}{_1}(a;b;z)$}
has integral representation
\begin{equation}
M(a,b,z) = \frac{1}{B(a,b-a)}\int_0^1 du\, e^{zu}u^{a-1}(1 - u)^{b - a- 1} \;.
\end{equation}
In terms of this function, we can write
\begin{equation}
\phi_1(s) = 1 - \frac{s}{s+1}{\cdot}\frac{M(s+1,s+2,\rho)}{M(s,s+1,\rho)}
= 1 - s{\cdot}\frac{M(1,2,\rho)}{M(0,1,\rho)}  + O(s^2) \;.
\end{equation}
Thus, the requirement that
\mbox{$\phi_1(0) = 1$}
is recovered and, noting that
\begin{equation}
M(1,2,\rho) = (e^\rho - 1)/\rho \;, \quad
     M(0,b,\rho) = 1 \;,
\end{equation}
we see that the mean BP for
\mbox{$N\to\infty$}
is
\mbox{$(e^\rho - 1)/\rho$}
as expected from the explicit exact result (\ref{BPMeanEx}).

Furthermore, a function that is entire in all arguments is obtained via
\begin{equation}
\mathbb{M}(a,b,z) \equiv M(a,b,z)/\Gamma(b) \;.
\end{equation}
With this definition and the relationship
\begin{equation}
\mathbb{M}(a,b,z) = e^z \mathbb{M}(b-a,b,-z) \;,
\end{equation}
we arrive at the result
\begin{equation}
\phi_1^{(\alpha)}(s) = (\alpha+1)\frac{\mathbb{M}(s,s+\alpha+2,\rho)}
     {\mathbb{M}(s,s+\alpha+1,\rho)} \;.
\label{MRatio}
\end{equation}
In this expression, both numerator and denominator are entire functions,
and we may consequently observe that the zeros of $\mathbb{M}(s,s+1,\rho)$
are identified with the eigenvalues of the full problem, while the zeros
of $\mathbb{M}(s,s+2,\rho)$
are identified with the eigenvalues of the truncated problem.
In fact, it follows directly from Cramer's rule that
\begin{equation}
\phi_1(s) = \langle e_1|A^{-1}(s)|e_1\rangle = \det A_{11}(s)/\det A(s) \;
\end{equation}
where $A_{11}$ denotes $A$ with the first row and column removed.

In view of the foregoing discussion, we may write the BP MGF as
\begin{equation}
\phi_1(s) \simeq \frac{1}{1 - s/\chi_1}\prod_{\ell = 2}^L
     \frac{1 - s/\chi^{(1)}_{\ell-1}}{1 - s/\chi_\ell} \;,
\end{equation}
where an exact expression is produced when
\mbox{$L = \infty$}.
This product form of linear factors is directly implied by (\ref{MRatio}),
but we also know that
the $\chi_\ell$, $\chi^{(1)}_\ell$ correspond to the eigenvalues of the matrices
$A$, $A_{11}$, respectively.
The two sequences of eigenvalues (ordered according to increasing modulus)
are interleaved, and both tend to
negative integral values as $\ell$ increases.
Moreover, they can be paired up in such a way that
\mbox{$\chi_{\ell+1} - \chi^{(1)}_\ell \to 0$}
as
\mbox{$\ell\to\infty$}.
The following inequalities hold:
\begin{equation}
0 < -\chi^{(1)}_\ell < -\chi_{\ell+1} < \ell \;,
\end{equation}
for
\mbox{$\ell = 1,2,\ldots$}.
Therefore, the individual ratios in the product above eventually become indistinguishable
from unity, implying that only a finite number, $L$, needs to be retained.
The residue theorem may be applied to recover the SF
as an exponential mixture
\begin{equation}
\bar{F}(t) \simeq \sum_{k = 1}^L W_k e^{\chi_k t} \;,
\end{equation}
with the weights being given by
\begin{equation}
W_1 = -\chi_1\prod_{\ell = 2}^L
     \frac{1 - \chi_1/\chi^{(1)}_{\ell-1}}{1 - \chi_1/\chi_\ell} \;, \quad
W_k = -\chi_k{\cdot}\frac{1 - \chi_k/\chi^{(1)}_{k-1}}{1 - \chi_1/\chi_k}
     \prod_{\ell = 2 \atop \ell\neq k}^L
     \frac{1 - \chi_k/\chi^{(1)}_{\ell-1}}{1 - \chi_k/\chi_\ell} \;,
\end{equation}
for
\mbox{$k = 2,3,\ldots$}.
The weights $W_k$ are positive, sum to unity and tend to zero for large $k$.
It follows that,
if we set
\mbox{$U_k \equiv -1/\chi_k$},
then the node/weight pairs
\mbox{$(U_k,W_k)$}
define a texture distribution of the sort discussed in Section~\ref{Spectral}.

\begin{figure}
\FIGURE
{\includegraphics[width=\wscl\linewidth, height=\hscl\linewidth]{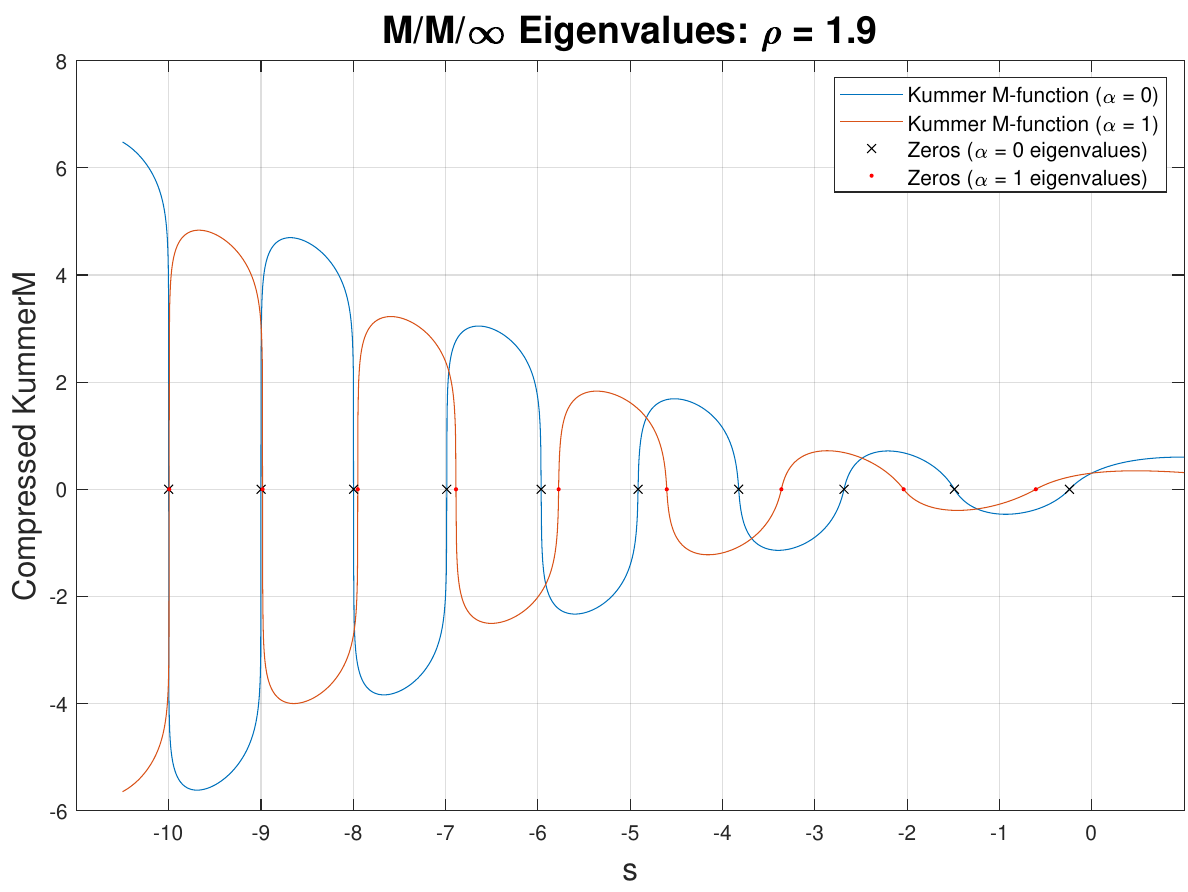}}
{\hphantom{x}\label{KummerZero}}
{Zeros of the Kummer M-functions in the numerator and denominator of the M/M/$\infty$
     BP MGF for partial traffic intensity $\rho = 1.9$.
     Denominator zeros are indicated by crosses. Numerator zeros are indicated by dots.}
\end{figure}

The $\chi_{\ell}$ are bracketed according to
\mbox{$\ell-1 < -\chi_\ell < \ell$}
for all
\mbox{$\ell = 1,2,\ldots$}
while, for the $\chi^{(1)}_{\ell}$, we may write
\begin{align}
\begin{aligned}
\ell-1 < -\chi^{(1)}_\ell \leq \ell   &\quad\text{for}\quad \ell = 1,\ldots,\floor(\rho) \;, \\
\ell   \leq -\chi^{(1)}_\ell < \ell+1 &\quad\text{for}\quad \ell = \ceil(\rho),\ldots \;, \\
-\chi^{(1)}_\rho = \rho               &\quad\text{for}\quad \text{integral $\rho \geq 1$} \;.
\end{aligned}
\end{align}
This behaviour can be observed in Figure~\ref{KummerZero} for the case
\mbox{$\rho = 1.9$}.
The function values on the vertical axis have been logarithmically compressed such the values
indicate the exponent of the order of magnitude.
The bijective compression function that was employed is given by
\begin{equation}
f(x) \equiv \sgn(x)\log_{10}(1 + |x|) = x{\cdot}\frac{\log_{10}(1 + |x|)}{|x|} \;,
\end{equation}
where $\sgn(x)$ denotes the sign function. Thus,
\mbox{$y = f(x)$}
can be inverted according to
\mbox{$x = \sgn(y)(10^{|y|} - 1)$}.
A bracketed root-finding method, such as bisection or a secant method, can be used to compute
the desired number of zeros of the relevant Kummer functions.
Alternatively, an eigenvalue problem can be solved for the matrices $A$, $A_{11}$ truncated to a
sufficiently large finite dimension (which will be greater than $L$).
An initial dimension can be estimated and subsequently iterated until convergence is achieved.

\begin{figure}
\FIGURE
{\includegraphics[width=\wscl\linewidth, height=\hscl\linewidth]{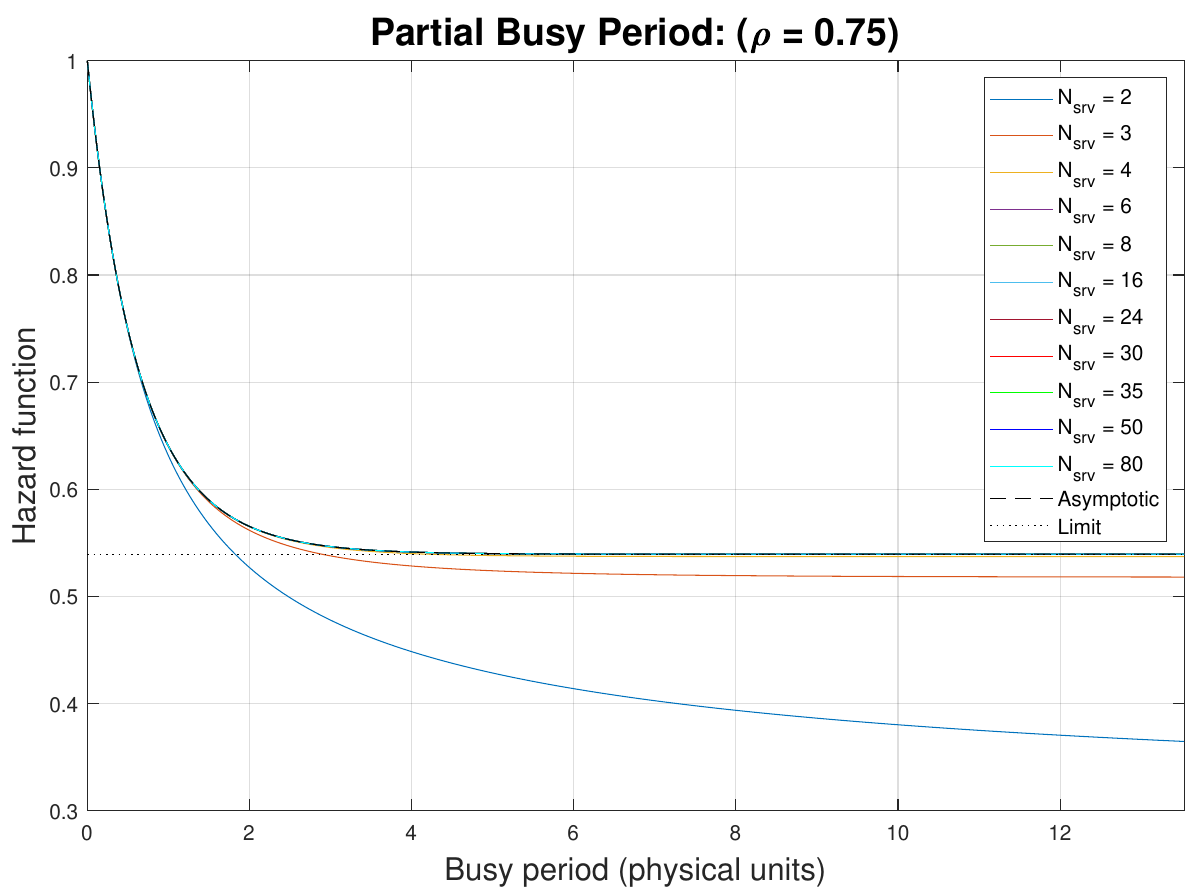}}
{\hphantom{x}\label{BPHazard_MultiN}}
{BP hazard function for various numbers of servers ($N_{\text{srv}}$) and
     partial traffic intensity $\rho = 0.75$. The dashed black curve is the hazard function
     for the M/M/$\infty$ model. The dotted black line is large-$t$ limit and lower bound of
     the asymptotic hazard function.}
\end{figure}

The eigenvalue of least modulus will dominate the large-t behaviour of the distribution
and provide a lower bound to the hazard function
\mbox{$H(t) \equiv P(t)/\bar{F}(t) > |\chi_1|$}.
This asymptotic level is indicated in Figure~\ref{BPHazard_MultiN} as the dotted line.
The hazard function for the M/M/$\infty$ model corresponds to the dashed black curve.
It is clear that this curve is approached rapidly as the number of servers increases.

\subsection{Constant-$\mathbf{r}$ Model}
We shall first consider the boundary case
\mbox{$r = 1$}
when
\mbox{$N \gg 1$}.
Let
\mbox{$m \equiv (\Lambda/\mu_1)^2$}
with $\Lambda$ given by (\ref{Lambda}),
noting that
\mbox{$m\to\infty$}
as
\mbox{$N\to\infty$}
since, by Stirling's formula,
\mbox{$m \asym{N \gg 1} e^{2N}/(2\pi N)$}.
We already know that
\mbox{$\bar{F}(mt) \asym{t\to\infty} 1/\sqrt{\pi t}$}
for any $N$, and we can write
\begin{equation}
\bar{F}_{\text{cut}}(mt) = \frac{2\mu_1}{\pi}\int_0^{4m} \frac{du}{4m}\, e^{-ut}
     R_{\text{cut}}(-u/m) \sqrt{\frac{4m - u}{u}} \;.
\end{equation}
Thus,
\begin{equation}
\bar{F}(mt) \asym{N \gg 1} \frac{\mu_1}{\pi}\int_0^\infty
     \frac{du}{\sqrt{u}}\, e^{-ut}\frac{1}{\sqrt{m}}R_{\text{cut}}(-u/m) \;.
\end{equation}
But also
\begin{equation}
\frac{\mu_1}{\sqrt{m}}R_{\text{cut}}(-u/m) \asym{u/m \ll 1} \frac{1}{1 + u} \;,
\end{equation}
which leads to the result
\begin{equation}
\bar{F}(mt) \asym{N \gg 1} \frac{1}{\sqrt{\pi}}\int_0^\infty \frac{dx}{\sqrt{x}}\,
     e^{-x}\frac{1}{\sqrt{1 + t/x}}
= \erfcx(\sqrt{t}) \;,
\end{equation}
where
\mbox{$\erfcx(t) \equiv e^{t^2}{\cdot}\erfc(t)$}
is the scaled complementary error function,
and $\erfc(t)$ denotes the standard complementary error function.
In completing this derivation, we have used the identity for the
incomplete gamma function
\mbox{$\Gamma(\half,t) = \sqrt{\pi}\erfc(\sqrt{t})$}.
The former integral expression is amenable to efficient numeral evaluation
via Gauss-Laguerre quadrature.
The asymptotic BP PDF is given by
\begin{equation}
P(t) \asym{N \gg 1} \frac{1}{\sqrt{\pi mt}} - \erfcx(\sqrt{t/m}) \;.
\end{equation}
The log-mean is given by
\begin{equation}
\langle\ln T_{\text{bp}}\rangle_{T_{\text{bp}}} \asym{N \gg 1} \ln m -\gamma_{\text{e}} \;.
\end{equation}
This result follows easily from direct integration of the representation
\begin{equation}
m{\cdot}P(mt) \asym{N\gg 1} \frac{2}{\pi}\int_0^\infty dv\, \frac{v^2}{1 + v^2} e^{-tv^2} \;.
\end{equation}

\begin{figure}
\FIGURE
{\includegraphics[width=\wscl\linewidth, height=\hscl\linewidth]{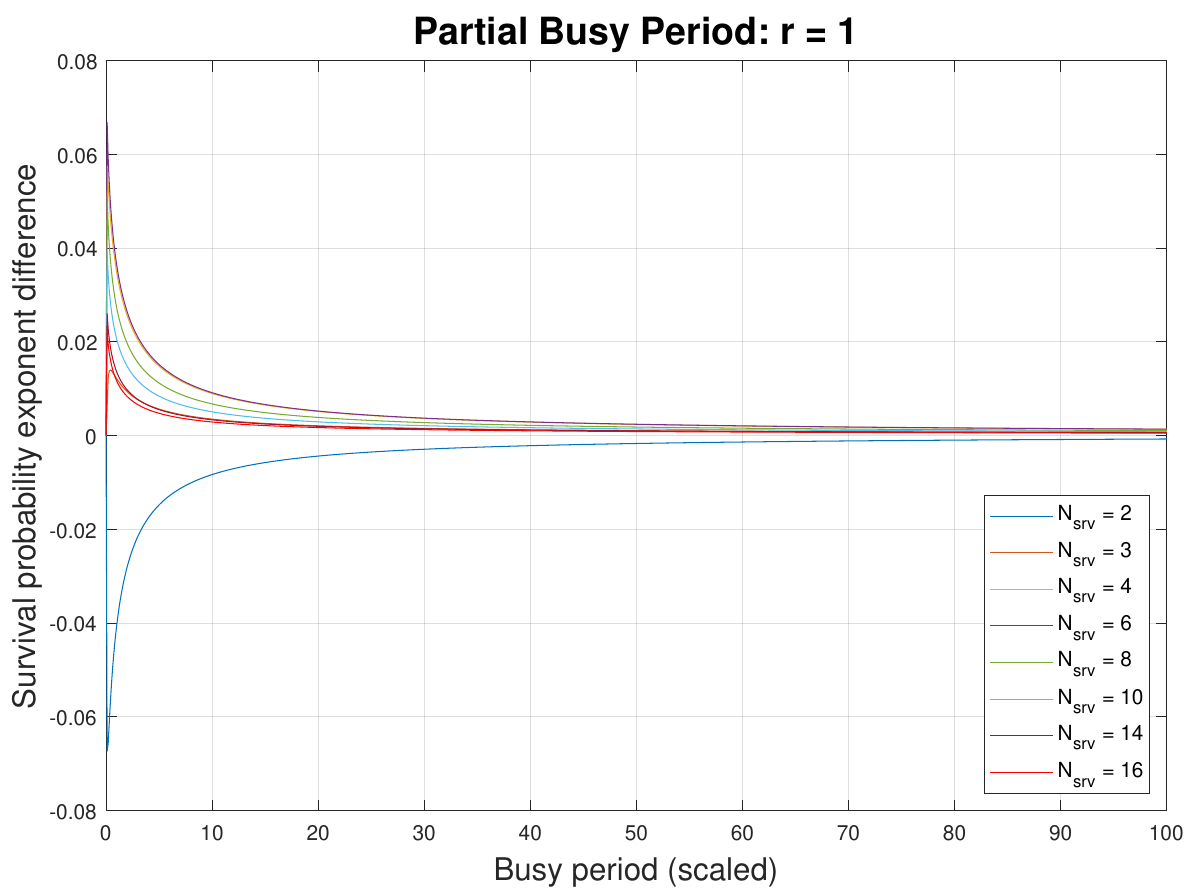}}
{\hphantom{x}\label{BPSFR1}}
{Deviation of the BP survival function from its large-$N$ asymptotic form for
     various numbers of servers ($N_{\text{srv}}$) and traffic intensity $r = 1$.
     The difference of the survival function logarithms is plotted on the vertical axis. Time on the horizontal axis
     is measured in units of the scale $m = (\Lambda/\mu_1)^2$.}
\end{figure}

In the ergodic region
\mbox{$0 < r < 1$},
the BP distribution tends to a mixture of two exponentials as $N$ become large:
one that dominates in the tail for large times and one that dominates for very small times.
For a fixed number of servers $N$, there is always a critical total traffic intensity
\mbox{$r = r_{\text{c}}$}
beyond which there is no contribution to the distribution from discrete poles.
On the other hand, for any fixed traffic intensity $r$, there exists at least one discrete pole
when the number of servers is sufficiently large, and the position of smallest-magnitude pole $\zeta_1$
tends to zero as the number of servers increases.
The exponential term generated by smallest-magnitude pole also dominates the large-$t$ behaviour
of the distribution.
Given that $\zeta_1$ approaches zero for
\mbox{$N \gg 1$},
it follows that this exponential asymptotic behaviour will set in ever earlier in $t$ as $N$
continues to increase. Hence, we must have that
\begin{equation}
\bar{F}(m_{\text{bp}}{\cdot}t) \asym{N\gg 1} \nu e^{-t}
\end{equation}
for all $t$ bounded away from zero, where
\mbox{$m_{\text{bp}} \equiv \langle T_{\text{bp}}\rangle_{T_{\text{bp}}}$}
denotes the mean BP.
A corollary of this is that
\mbox{$\zeta_1 \asym{N\gg 1} -1/m_{\text{bp}}$},
which is easily verified numerically by computing $\zeta_1$ from (\ref{PoleNearZero})
and comparing with the exact result for $m_{\text{bp}}$.
The constant $\nu$ is derived from the residue of the pole at $s = \zeta_1$, namely
\mbox{$\nu = -\mu_1 J(\zeta_1)/\zeta_1$},
where the function $J(s)$ has been given in (\ref{JFn}).
It also holds that
\mbox{$\nu \to 1^-$}
as
\mbox{$N\to\infty$}.
On the timescale that corresponds to measuring $t$ in units of the mean BP, we have
a PDF of the mixture form
\begin{equation}
m_{\text{bp}} P(m_{\text{bp}}{\cdot}t) \asym{N\gg 1} (1 - \nu)\delta(t) + \nu e^{-t} \;,
\end{equation}
where
\mbox{$\delta(t)$}
denotes the Dirac delta function.
The SF for the two-exponential mixture that leads to this PDF is given by
\begin{equation}
\bar{F}(t) \asym{N\gg 1} (1 - \nu) e^{-t/m''} + \nu e^{-\nu t/m'} \;,
\label{SFAsymN}
\end{equation}
with
\mbox{$m' = m_{\text{bp}} - (1-\nu)m''$}
chosen to reproduce the correct mean, and $m''$ is a constant of order unity
estimated as
\mbox{$m'' \simeq (1-\nu)/P(0) = (1 - \nu)N$}.
It follows that
\mbox{$\nu/m' \asym{N\gg 1} -\zeta_1$},
which strengthens the relationship given above.
It is easily confirmed numerically that this result remains accurate
down to quite small values of $N$.

\begin{figure}
\FIGURE
{\includegraphics[width=\wscl\linewidth, height=\hscl\linewidth]{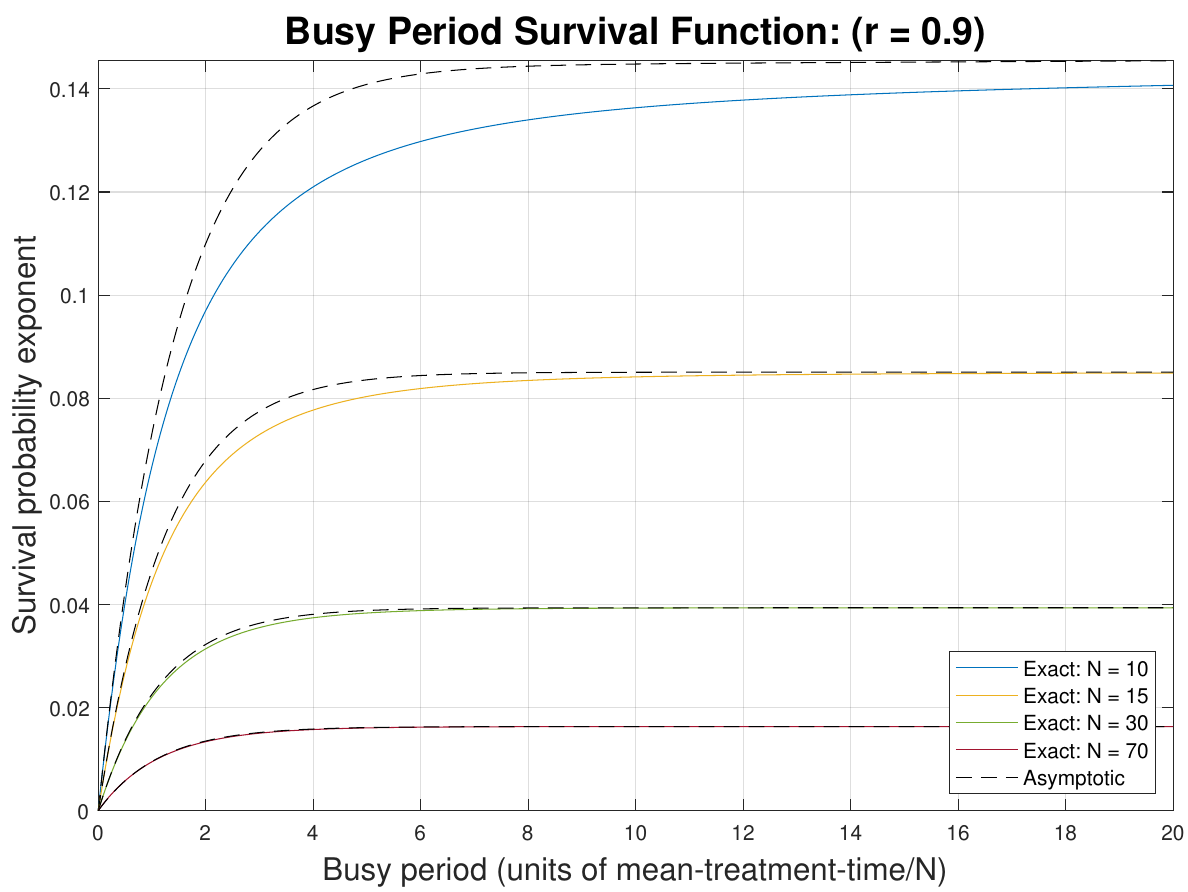}}
{\hphantom{x}\label{BPAsymTest}}
{BP survival functions for various numbers of servers ($N$) and
     traffic intensity $r = 0.9$, compared with their respective large-$N$ asymptotic forms.
     These curves elucidate the very short-time behaviour. The asymptotic forms are evaluated
     at the corresponding values of $N$.}
\end{figure}

In Figure~\ref{BPAsymTest}, we have plotted the SF exponents
$-\log_{10}\bar{F}(t)$ for an increasing sequence of server numbers $N$,
where values on the time axis correspond directly to
the model definition in the introductory section with the choice
\mbox{$\mu = 1/N$}
adopted there.
Thus, time is being measured in units of $1/N$ times the mean treatment time.
If we were to consider a situation where both $r$ and $\mu$ were to remain constant with respect
to $N$, in which case $\lambda$ would be linear in $N$, then
any given value on the time axis would be a different point in physical time for each curve.
It can be seen from the figure that, in the short-time regime being plotted,
agreement between the numerical exact curves and the
analytical asymptotic curves  improves rapidly with increasing $N$.

\begin{figure}
\FIGURE
{\includegraphics[width=\wscl\linewidth, height=\hscl\linewidth]{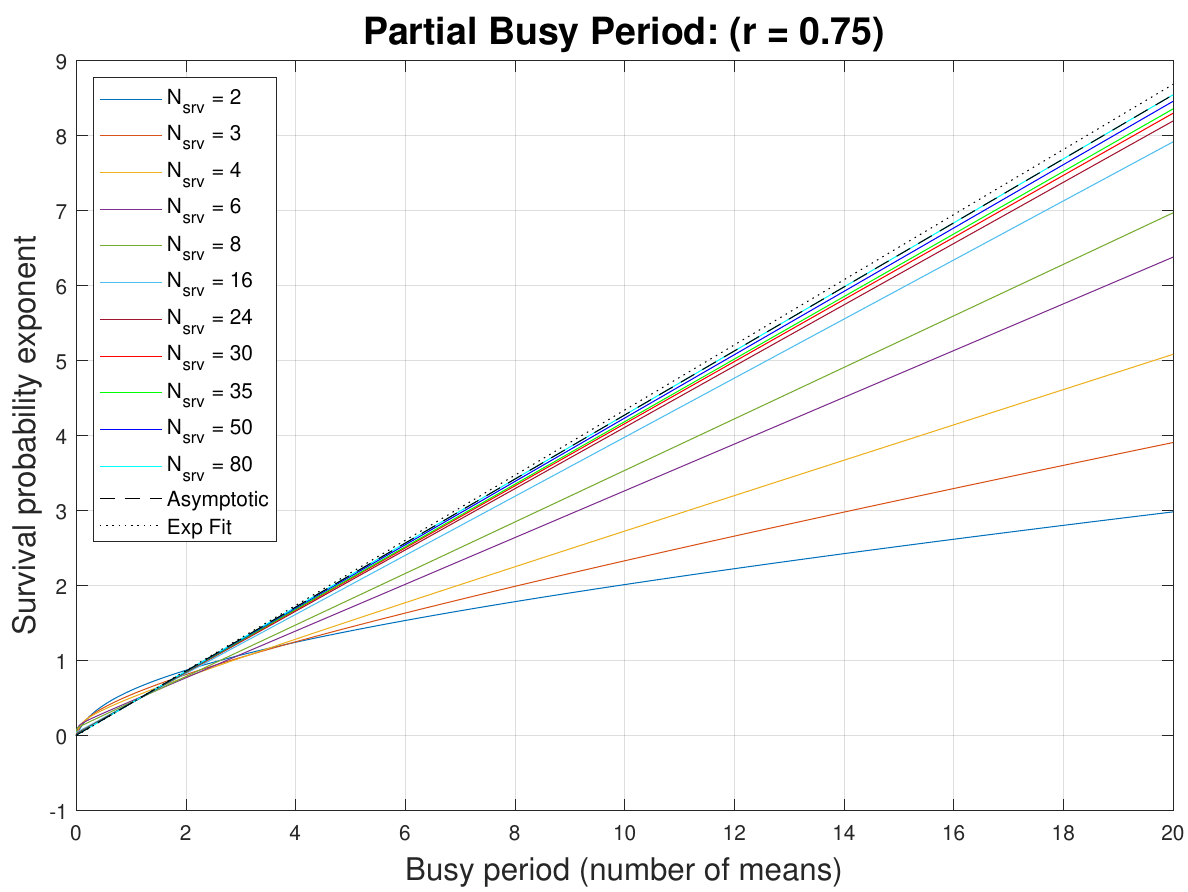}}
{\hphantom{x}\label{BPSFTail_MultiNRConst}}
{BP survival functions for various numbers of servers ($N_{\text{srv}}$) and
     traffic intensity $r = 0.75$, showing long-time behaviour. The dashed black line is the asymptotic form
     evaluated for $N = 80$. The black dotted line is the strict $N\to\infty$ limiting exponential.}
\end{figure}

The long-time behaviour is presented in Figure~\ref{BPSFTail_MultiNRConst},
where values on the time axis are measured in units of the mean BP.
The black dashed curve is the analytical asymptotic curve evaluated here for
\mbox{$N = 80$},
and is seen to coincide with its numerical exact counterpart to within the
line-width of the graph.
The black dotted line is the limiting exponential that represents the strict
\mbox{$N\to\infty$}
limit
\mbox{$\bar{F}(m_{\text{bp}}{\cdot}t) \asym{N\to\infty} e^{-t}$}.
It constitutes a lower bound to the tail of the SFs for all $N$,
and is approached from above as
\mbox{$N\to\infty$}.

\section{Complex-Pole Method}
\label{Complex}
The explicit expression for the cut function developed in the previous section
allows us to bypass the spectral decomposition and derive explicit closed-form
results for the cut contribution to the BP PDF for an arbitrary number of servers.
These comprise sums of Bessel and Marcum Q-functions
where the individual terms are generated by
poles in a complex plane that are induced by the
zeros of the cut polynomial.
In fact, it will become clear that knowledge of the cut polynomial
is all that is required to completely determine the BP distribution
for any given number of servers.
While this method is less robust numerically than the spectral approach, it
provides more insight into the analytic structure of the problem.

The general form of the cut contribution to the BP PDF van be expressed as
\begin{equation}
P_{\text{cut}}(t) = -\frac{r\mu_1}{2}e^{-(1+r)t}\oint_\mathcal{C}\frac{dz}{2\pi iz}\,
     \left(z - \frac{1}{z}\right)^2e^{\sqrt{r}t(z + 1/z)}{\cdot}
     R_{\text{cut}}(s) \;,
\label{PcutZ}
\end{equation}
where the contour $\mathcal{C}$
denotes the anti-clockwise traversed unit circle about the origin, and
\mbox{$R_{\text{cut}}(s)$}
is the cut function, with
\begin{equation}
s      = \sigma - r - \mu_1 \;, \quad
\sigma = \sqrt{r}(z + 1/z) - 1 + \mu_1 \;.
\end{equation}
This follows from setting
\mbox{$v = \cos\theta$}
in (\ref{FcutV}), and then writing
\mbox{$\cos\theta = -(z + 1/z)/2$}
where $z$ lies on the unit circle in the complex
$z$-plane\footnote{The minus has been included in order to reproduce the
standard sign in the exponential term of the Schaefli formula for the modified
Bessel function ({\it c.f.}\ (\ref{Schlaefli})).}.
For
\mbox{$N = 1$}
(\mbox{$\mu_1 = 1$}),
the cut function is the trivial constant
\mbox{$R_{\text{cut}}(s) = 1/r$}.
For
\mbox{$N = 2$}
(\mbox{$\mu_1 = 1/2$}),
\begin{equation}
R_{\text{cut}}(s) = \frac{1}{r - \mu_1\sigma} \;,
\end{equation}
while, for
\mbox{$N = 3$}
(\mbox{$\mu_1 = 1/3$}),
\begin{equation}
R_{\text{cut}}(s) = -\frac{2r}{\sigma^3 + (\mu_1 - r)\sigma^2 - 4r^2\mu_1} \;.
\end{equation}
By writing the cut polynomial in terms of its roots,
\mbox{$C_N(\sigma) = \prod_{\ell=1}^{2N-3} (\sigma - \sigma^{(\ell)})$},
we have the general formula
\begin{equation}
R_{\text{cut}}(s) = -\frac{1}{r\mu_1^2}{\cdot}\prod_{k=1}^{N-1}(r\mu_k)
     {\cdot}\frac{1}{\prod_{\ell=1}^{2N-3} (\sigma - \sigma^{(\ell)})} \;.
\end{equation}

The Schlaefli formula for the modified Bessel function of the first kind is given by
\begin{equation}
I_n(t) = \oint\frac{dz}{2\pi i}\, \frac{1}{z^{n+1}} \exp\left\{\frac{t}{2}\left(z + \frac{1}{z}\right)\right\} \;,
\label{Schlaefli}
\end{equation}
where the closed contour is taken around the origin.
This may be used to derive the generating function
\begin{equation}
\sum_{n=-\infty}^{+\infty} z^nI_n(t) = \exp\left\{\frac{t}{2}\left(z + \frac{1}{z}\right)\right\} \;,
\end{equation}
from which it follows that
\begin{equation}
\exp\left\{\frac{t}{2}\left(z + \frac{1}{z}\right)\right\} = I_0(t) + \sum_{n=1}^{\infty}
     \left(z^n + z^{-n}\right)I_n(t) \;.
\end{equation}
Thus, suppose that
\mbox{$g(z)$}
is a meromorphic function such that
\mbox{$g(z) = g(1/z)$}
and that $\mathcal{C}$ denotes the anti-clockwise unit circle. Then we have the integral
\begin{align}
\begin{aligned}
\mathcal{I} &\equiv \oint_\mathcal{C}\frac{dz}{2\pi iz}\, g(z)
     \exp\left\{\frac{t}{2}\left(z + \frac{1}{z}\right)\right\} \\
&= \oint_\mathcal{C}\frac{dz}{2\pi iz}\, g(z) {\cdot} I_0(t) +
     2\sum_{n=1}^{\infty}I_n(t)\oint_\mathcal{C}\frac{dz}{2\pi iz}\, g(z)z^n \;.
\end{aligned}
\label{IDef}
\end{align}
Suppose now that
\mbox{$G(z) \equiv g(z)/z$}
has only simple poles within the unit circle at
\mbox{$z = \beta_1,\beta_2,\ldots,\beta_M$}.
Then we have
\begin{equation}
\mathcal{I} = c_0 {\cdot}I_0(t) + 2\sum_{n=1}^{\infty} c_n I_n(t) \;,
\label{IIntegral}
\end{equation}
where the coefficients, for
\mbox{$n = 0,1,\ldots$},
are the residue sums
\begin{equation}
c_n \equiv \oint_\mathcal{C}\frac{dz}{2\pi i}\, z^n G(z) = \sum_{m=1}^M \beta_m^n \Residue[G(z),\beta_m] \;.
\label{cn}
\end{equation}

For
\mbox{$N = 1$},
(\ref{PcutZ}) and the Schaefli formula (\ref{Schlaefli}) reproduce the
well-known result \citep{BP:Artalejo01}
\begin{equation}
P_{\text{cut}}(t) = e^{-(1+r)t}\left[I_0(2\sqrt{r}t) - I_2(2\sqrt{r}t)\right]
= \frac{1}{\sqrt{r}t}e^{-(1+r)t}I_1(2\sqrt{r}t) \;.
\end{equation}
In the general case, for
\mbox{$N > 1$},
the relationships
\begin{equation}
C_N(\sigma) = \prod_{\ell=1}^{2N-3} \left(\sigma - \sigma^{(\ell)}\right) \;, \quad
\sigma      = \sqrt{r}(z + 1/z) +\mu_1 - 1 \;,
\end{equation}
lead to the form
\begin{equation}
C_N(\sigma) = \left(\frac{\sqrt{r}}{z}\right)^{2N-3}\prod_{\ell=1}^{2N-3}
     \left(z^2 + 2\alpha_\ell z + 1\right) \;,
\end{equation}
where
\mbox{$\alpha_\ell \equiv (\mu_1 - 1 - \sigma^{(\ell)})/(2\sqrt{r})$}.
Hence, (\ref{PcutZ}) becomes
\begin{equation}
P_{\text{cut}}(t) = \frac{\sqrt{r}}{2}\left({\textstyle \prod_{k=2}^{N-1}\mu_k}\right)
     e^{-(1+r)t}\oint\frac{dz}{2\pi i}\, e^{\sqrt{r}t(z + 1/z)} G(z) \;,
\end{equation}
where (\ref{IDef}) applies with
\begin{equation}
G(z) = \frac{g(z)}{z} = z^{2(N-3)}{\cdot}
     \frac{\left(1 - z^2\right)^2}{\prod_{\ell=1}^{2N-3}\left(z^2 + 2\alpha_\ell z + 1\right)} \;,
\end{equation}
which allows us to the identify the $c_n$ in the relevant instance of (\ref{IIntegral}).
We note that for
\mbox{$N \geq 3$}
there is no pole at the origin.

Next, we write
\begin{equation}
z^2 + 2\alpha_\ell z + 1 = (z - \beta_\ell)(z - 1/\beta_\ell) \;,
\end{equation}
where $\beta_\ell$ is chosen to be the unique root of the pair of roots
\mbox{$\beta_\ell = -\alpha_\ell \pm\sqrt{\alpha_\ell^2 - 1}$}
that lies within the unit circle
\mbox{$|\beta_\ell| < 1$}.
In each pair, one always lies inside and one outside the unit circle,
and they are complex conjugates.
With this in mind, we can express (\ref{cn}) as
\begin{equation}
c_n = -\sum_{k=1}^{2N-3} \beta_k^n\gamma_k \;, \quad
     \gamma_k \equiv \beta_k^{2N-5}(1 - \beta_k^2)/D_k \;,
\end{equation}
with
\begin{equation}
D_k \equiv \prod_{\ell = 1 \atop \ell \neq k}^{2N-3} (\beta_k - \beta_\ell)
     (\beta_k - 1/\beta_\ell) \;.
\end{equation}
In terms of these quantities, the cut PDF for
\mbox{$N \geq 3$}
reads
\begin{equation}
P_{\text{cut}}(t) = \frac{\sqrt{r}}{2}\left({\textstyle \prod_{k=2}^{N-1}\mu_k}\right)
     e^{-(1+r)t}\left[
     -c_0{\cdot}I_0(\sqrt{4r}t) + 2\sum_{k=1}^{2N-3}\gamma_k{\cdot}
     \sum_{n=1}^{\infty} \beta_k^n I_n(\sqrt{4r}t)\right] \;.
\label{PcutI}
\end{equation}
It is useful to recall that
\mbox{$c_0 = -\sum_{k=1}^{2N-3}\gamma_k$}.

The generalized Marcum Q-function for order $\nu$ is defined by the integral
\begin{equation}
Q_\nu(a,b) \equiv \frac{1}{a^{\nu-1}}\int_0^b dx\, x^\nu \exp\left\{\half(x^2 + a^2)\right\}
     I_{\nu-1}(x) \;,
\end{equation}
where $I_\nu(x)$ is the modified Bessel function of the first kind of order $\nu$.
It constitutes a cumulative distribution function (CDF) in the variable $b$.
The original Marcum Q-function
\mbox{$Q(a,b)$}
is the special case of the generalized Q-function
\mbox{$Q_\nu(a,b)$}
for
\mbox{$\nu = 1$}.
Its complementary function
\mbox{$\bar{Q}(a,b) \equiv 1 - Q(a,b)$}
can be represented by the infinite Neumann expansion
\begin{equation}
\bar{Q}(a,b) = e^{-(a^2+b^2)/2}\sum_{\alpha = 1}^\infty \left(\frac{b}{a}\right)^\alpha
     I_\alpha(ab) \;,
\label{Neu1}
\end{equation}
which allows it to be extended to complex-valued arguments $a, b$.
The foregoing series is useful for numerical computation when
\mbox{$|b/a| < 1$}.
Otherwise, one may appeal to the alternative form
\begin{equation}
\bar{Q}(a,b) = 1 - e^{-(a^2+b^2)/2}\sum_{\alpha = 0}^\infty \left(\frac{a}{b}\right)^\alpha
     I_\alpha(ab) \;.
\label{Neu2}
\end{equation}
The Marcum Q-function was originally introduced in radar detection theory for
non-fluctuating targets \citep{BP:Marcum60},
and has subsequently found applications in communications and signal processing.
\citet{BP:Abate89} were first to use it in the context of queueing theory.
Here we adopt it
to write (\ref{PcutI}) in the more compact form
\begin{equation}
P_{\text{cut}}(t) = \frac{\sqrt{r}}{2}\left({\textstyle \prod_{k=2}^{N-1}\mu_k}\right)
     \left[-c_0{\cdot}e^{-(r+1)t}I_0(\sqrt{4r}t) +
     2e^{-(r+\mu_1)t}\sum_{k=1}^{2N-3}\gamma_k
     {\cdot}e^{\sigma_k t}\bar{Q}\left(a_k(t),b_k(t)\right)\right]
\end{equation}
where
\begin{equation}
a_k(t) = \left(\sqrt{4r}t/\beta_k\right)^{1/2} \;, \quad
     b_k(t) = \left(\sqrt{4r}\beta_k t\right)^{1/2} \;,
\end{equation}
and these are in general complex valued.
A robust numerical algorithm for the computation of the Marcum
Q-function\footnote{Direct implementation of (\ref{Neu1}), (\ref{Neu2}) also yields goods results.}
that allows for complex arguments has been given in \citep{BP:Gil14}.

For
\mbox{$N = 2$},
there is an additional pole at the origin.
In this case, we have
\begin{equation}
c_n = \oint_\mathcal{C}\frac{dz}{2\pi i}\, z^{n-2}\frac{\left(1 - z^2\right)^2}
     {z^2 + 2\alpha_1 z + 1} \;,
\end{equation}
from which it is evident that a pole at the origin contributes when
\mbox{$n = 0,1$}.
We find that
\mbox{$c_0 = 2\beta_1$},
\mbox{$c_1 = \beta_1^2$}
and, for
\mbox{$n \geq 2$},
\mbox{$c_n = -\beta_1^{n-1}(1 - \beta_1^2)$}.
It follows that
\begin{equation}
P_{\text{cut}}(t) = \sqrt{r}e^{-(r+1)t}\left[\beta_1 I_0(\sqrt{4r}t)
     +  I_1(\sqrt{4r}t) - \left(\frac{1}{\beta_1} - \beta_1\right)
     \sum_{n=1}^\infty \beta_1^n I_n(\sqrt{4r}t)\right] \;,
\end{equation}
and it turns out that
\mbox{$\beta_1 = \min(\sqrt{4r},1/\sqrt{4r})$},
so that
\begin{equation}
\frac{1}{\beta_1} - \beta_1 = \left|\sqrt{4r} - \frac{1}{\sqrt{4r}}\right| \;.
\end{equation}
This is similar to the expression for the BP distribution as an infinite series
of Bessel functions that was previously derived by
\citet{BP:Arora64} for the two-server case.
In terms of the complementary Marcum Q-function,
the two-server result becomes
\begin{equation}
P_{\text{cut}}(t) = 2e^{-(r+1)t}\left[\min\left(r,\tfrac{1}{4}\right){\cdot}I_0(\sqrt{4r}t)
     +  I_1(\sqrt{4r}t) - \left|r - \tfrac{1}{4}\right|
     {\cdot}e^{2(r+1/4)t}\bar{Q}\left(a_1(t),b_1(t)\right)\right] \;,
\end{equation}
and we may note that
\begin{equation}
a_1(t) = \left\{
\begin{array}{ccc}
\sqrt{t}   & \quad \text{for} \quad & r \leq 1/4 \;, \\
\sqrt{4rt} & \quad \text{for} \quad & r > 1/4 \;,
\end{array}\right. \quad
b_1(t) = \left\{
\begin{array}{ccc}
\sqrt{4rt} & \quad \text{for} \quad & r \leq 1/4 \;, \\
\sqrt{t}   & \quad \text{for} \quad & r > 1/4 \;.
\end{array}\right.
\end{equation}
In this instance, the arguments of the Marcum Q-function are real.

\begin{figure}
\FIGURE
{\includegraphics[width=\wscl\linewidth, height=\hscl\linewidth]{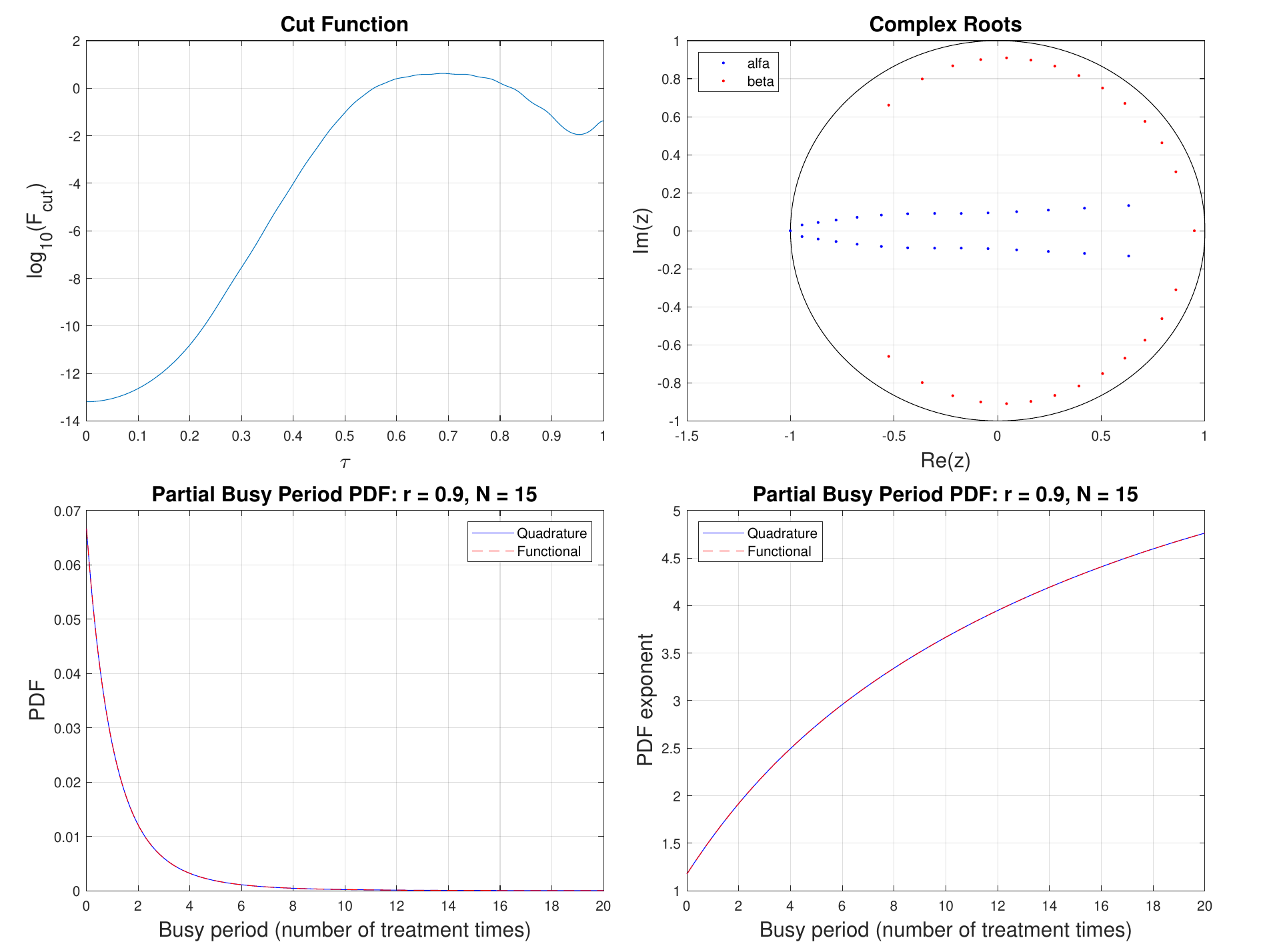}}
{\hphantom{x}\label{BPPDFNX}}
{BP distribution via the complex-pole method for $r = 0.9$, $N = 15$.
    The NW graph is the cut function $R_{\text{cut}}(\cdot)$ on a logarithmic scale as a function of $\tau$.
    The NE graph is the unit circle in the complex $z$-plane displaying the locations of the
    $\alpha$-poles (blue) and $\beta$-poles (red) within the unit circle.
    The single real $\alpha$-pole lies just outside the unit circle.
    The SW graph compares the survival function from the spectral method with that from the complex-pole method.
    The SE graph is the same comparison for the negative base-10 logarithm of the survival function and emphasizes agreement in the tail.}
\end{figure}

\section{Summary Statistics}
\label{Statistics}
\subsection{Analytical}
\begin{figure}
\FIGURE
{\includegraphics[width=\wscl\linewidth, height=\hscl\linewidth]{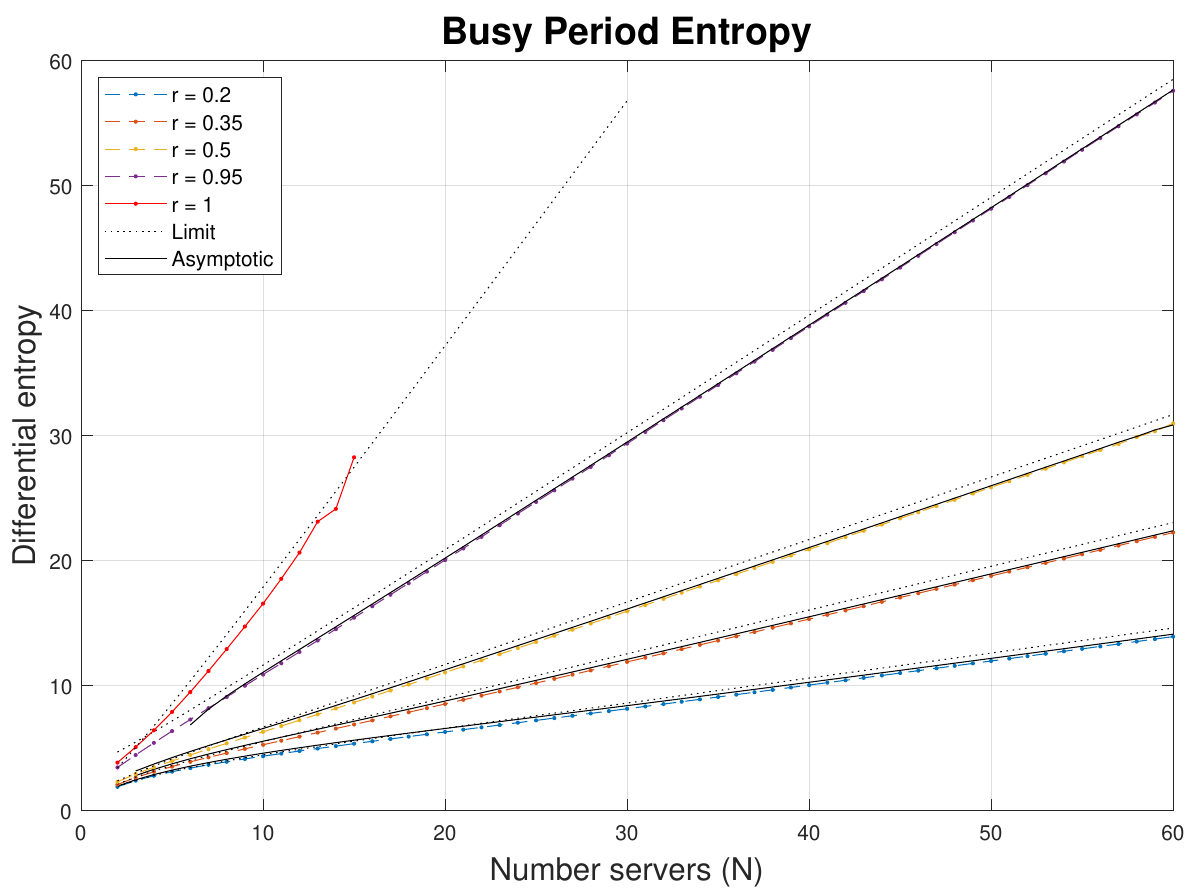}}
{\hphantom{x}\label{BPEntropy}}
{BP entropy versus number of servers ($N$) for various values of the traffic intensity (r).
     Each computed coloured dashed curve overlays a black solid curve that represents the analytical asymptotic form
     evaluated with the corresponding parameters $(r,N)$. The dotted black lines are the strict $N\to\infty$
     limiting curves for each value of $r$.}
\end{figure}

\begin{figure}
\FIGURE
{\includegraphics[width=\wscl\linewidth, height=\hscl\linewidth]{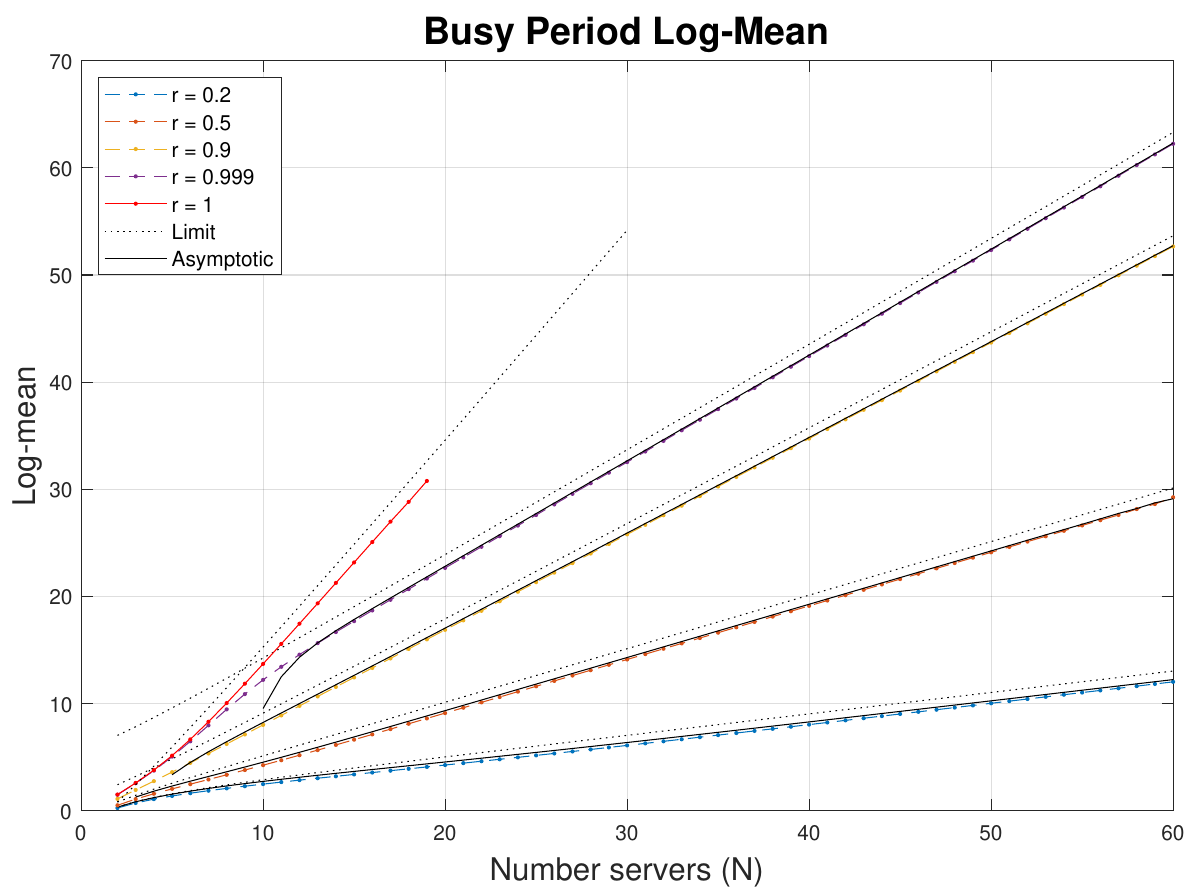}}
{\hphantom{x}\label{BPLogMean}}
{BP log-mean versus number of servers ($N$) for various values of the traffic intensity (r).
     Each computed coloured dashed curve overlays a black solid curve that represents the analytical asymptotic form
     evaluated with the corresponding parameters $(r,N)$. The dotted black lines are the strict $N\to\infty$
     limiting curves for each value of $r$.}
\end{figure}

\begin{figure}
\FIGURE
{\includegraphics[width=\wscl\linewidth, height=\hscl\linewidth]{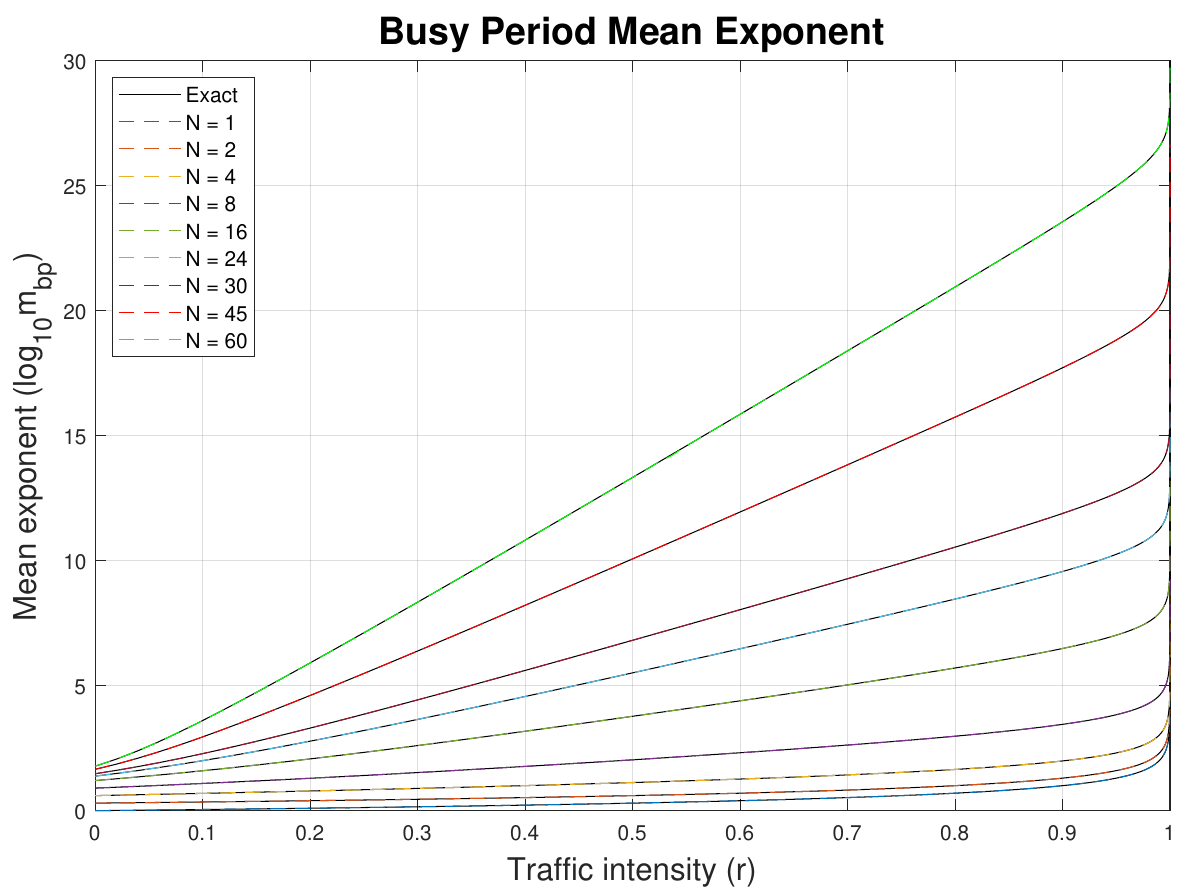}}
{\hphantom{x}\label{BPMeanExp}}
{BP mean versus traffic intensity $r$, plotted as its base-10 logarithm.
     For each number of servers $N$, the computed coloured dashed curve overlays a solid
     black curve representing the exact analytical result. Discrepancies are less than the linewidth.}
\end{figure}

In order to test how the numerical algorithms for computing the BP distribution,
developed in this work,
perform across a wide range of input parameters $(r, N)$ in the ergodic region,
we plot various summary statistics against a grid of these parameters.
We can then ascertain (i) whether the observed variations are sufficiently regular and continuous,
and (ii) how well the values approach the analytically calculated asymptotic limits
as the number of servers grows.
The results presented here are also useful in determining the values of $N$ beyond which the
asymptotic limit can serve as a proxy for the exact distribution, given a desired level of accuracy.

Since the mean and higher moment diverge rapidly as the traffic intensity $r$ approaches unity,
where heavy-tailed behaviour sets in,
we seek summary statistics that remain finite at
\mbox{$r= 1$}
and can be computed conveniently.
For this purpose, we have chosen the differential entropy
\mbox{$H \equiv \langle -\ln(P(T))\rangle_T$}
and the log-mean
\mbox{$L \equiv \langle\ln(T)\rangle_T$}.
The log-mean is widely used in radar detection theory to estimate the shape parameters
of candidate heavy-tailed distributions for high-resolution radar clutter from collected
experimental data \citep{BP:Rosen22}.

Figure~\ref{BPEntropy} presents the BP differential entropy plotted against number of servers $N$ for
various values of the traffic intensity $r$. It includes a comparison with asymptotic limits as detailed
in the caption.
Figure~\ref{BPLogMean} presents the BP log-mean plotted against number of servers $N$ for
various values of the traffic intensity $r$. It also includes a comparison with asymptotic limits as detailed
in the caption.
In Figure~\ref{BPMeanExp}, the BP mean is plotted against traffic intensity $r$ as its base-10 logarithm, for various
numbers of servers $N$, and is compared with the known exact results.
All these graphs demonstrate consistent behaviour for the numerical computations, agreement with exact results
and the expected approach to asymptotic limits.

\subsection{Simulation}
\begin{figure}
\FIGURE
{\includegraphics[width=\wscl\linewidth, height=\hscl\linewidth]{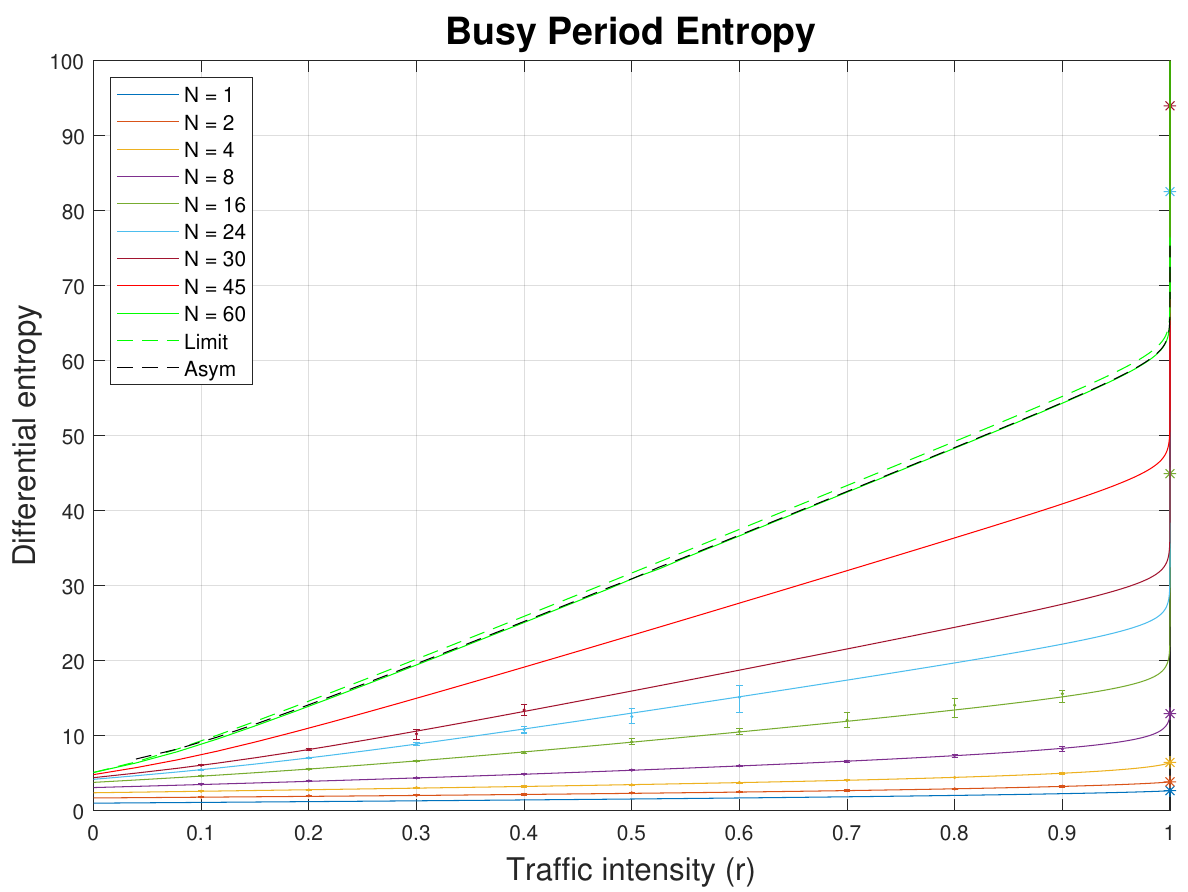}}
{\hphantom{x}\label{BPEntropySim}}
{BP entropy versus traffic intensity $r$ for various server numbers $N$,
     with overlaid simulation data as indicated by the dots with error bars.
     The green dashed curve is the extreme large-$N$ limit based on a single exponential and computed for $N = 60$ .
     The dashed black curve is the finer large-$N$ asymptotic form based on the two-exponential mixture and computed for $N = 60$.
     It aligns very closely with the exactly computed result.}
\end{figure}

\begin{figure}
\FIGURE
{\includegraphics[width=\wscl\linewidth, height=\hscl\linewidth]{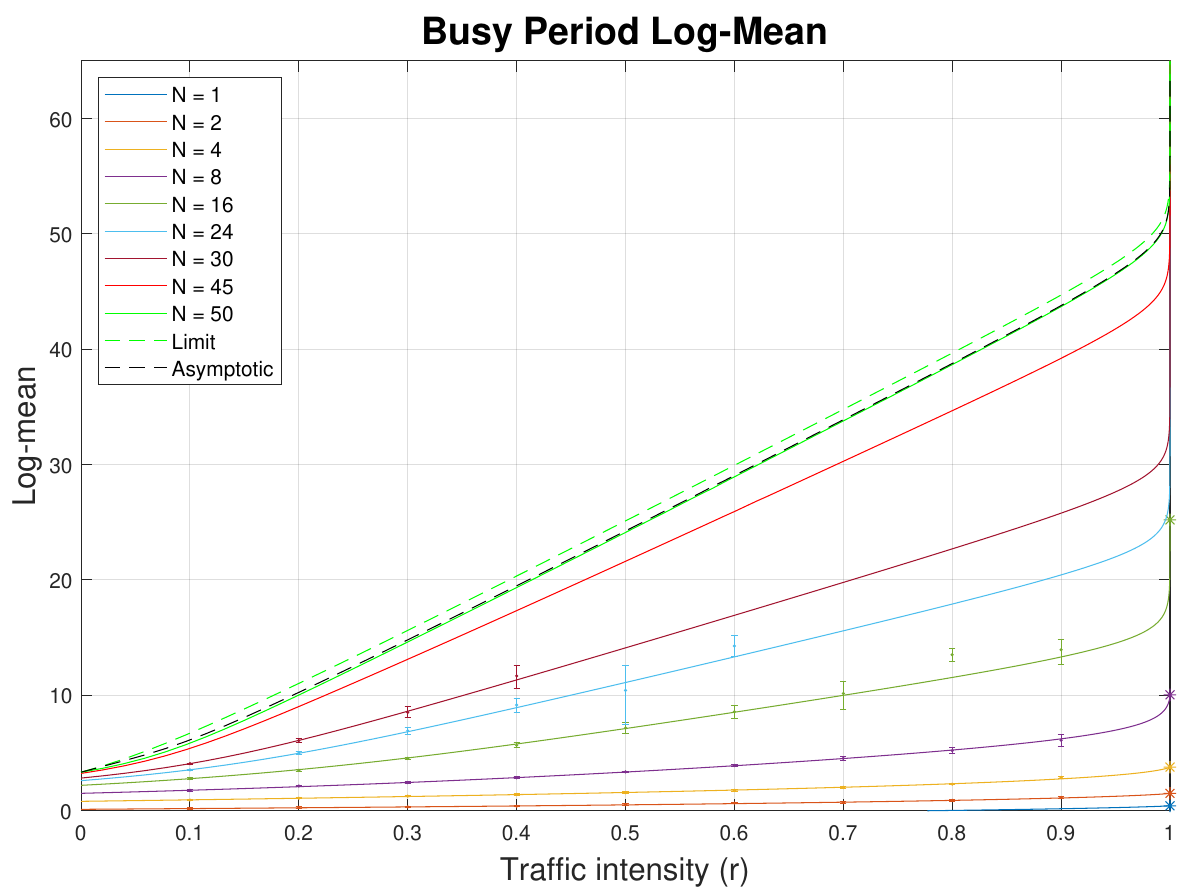}}
{\hphantom{x}\label{BPLogMeanSim}}
{BP log-mean versus traffic intensity $r$ for various server numbers $N$,
     with overlaid simulation data as indicated by the dots with error bars.
     The green dashed curve is the extreme large-$N$ limit based on a single exponential and computed for $N = 50$ .
     The dashed black curve is the finer large-$N$ asymptotic form based on the two-exponential mixture and computed for $N = 50$.
     It aligns very closely with the exactly computed result.}
\end{figure}

We have implemented a discrete event simulation (DES) for the M/M/$c$ system.
Adhering to the ergodic region
(\mbox{$0 < r < 1$}),
we run the simulation for a long time period as a steady-state simulation that is initialized
to the empty state, and collect the times (epochs) at which an empty system becomes non-empty
and at which a non-empty system becomes empty. These data are used to generate the
empirical partial BP distribution from which summary statistics can be estimated.
The goal is to compare the analytical results for various summary statistics with
their empirical estimates across a wide range of the input-parameter pairs
\mbox{$(r, N)$}
in order to verify the theoretical results derived in this work.
A simulation is run for a given parameter pair whenever the expected number of
BPs generated in the maximum tolerated simulation time $T_{\text{stop}}$
is at least $10$.
The maximum simulation time was taken to be $10^6$ multiples of the mean treatment time.
According to the theory of regenerative processes \citep{BP:Crane74}
the expected number of regeneration cycles (empty periods followed by a BP)
in time $T_{\text{stop}}$ is simply given by
\mbox{$n_{\text{reg}} = T_{\text{stop}}/[(m_{\text{bp}} + 1/r){\cdot}T_{\text{scl}}]$}
where $T_{\text{scl}}$ is a timescale needed to ensure that time in both the simulation
and the analysis are measured in the same units.
In the simulation, we set the mean treatment time to unity so that the simulation clock
runs in units of the mean treatment time.
In this case we must set
\mbox{$T_{\text{scl}} = 1/N$}.

The nearest-neighbour estimator \citep{BP:Kozachenko87,BP:Tsybakov96}
is used for the empirical differential entropy.
This choice avoids calculation of a kernel density estimator (KDE) for the
empirical PDF that is required by other empirical entropy estimators \citep{BP:Beirlant01}.
The KDE is problematic for distributions with semi-infinite support because
undesirable artefacts are inevitably produced at the boundary.
For a sample of size $M$
\mbox{$\{T_i: 1 \leq i \leq M\}$},
the nearest-neighbour estimator is given by
\begin{equation}
\hat{H} \equiv \frac{1}{M}\sum_{i=1}^M \ln[(M-1)\rho_i] + \ln 2 + \gamma_{\text{e}} \;,
\end{equation}
where
\mbox{$\rho_i \equiv \min_{j\neq i}\left\| T_i - T_j\right\|$}
is the nearest-neighbour Euclidean distance of $T_i$ from all other members of the sample.

Due to the small sample numbers in many cases ({\it viz.}\ $r$ close to unity, or large $N$)
and lack of alternatives for the entropy, the bootstrap method \citep{BP:Hesterberg11}
has been used to generate the
confidence intervals that provide the error bars in the graphs.
The significance level is taken to be
\mbox{$\alpha = 0.01$}.
Thus, the displayed error bars in Figures \ref{BPEntropySim} and \ref{BPLogMeanSim},
for the differential entropy and log-mean, respectively,
indicate the central $99\%$ confidence intervals.
In cases where only a dot appears, the length of the confidence interval is less than the dot size.

The nearest-neighbour entropy estimator is somewhat problematic for the bootstrap method
as it generates duplicates, for which the nearest-neighbour distance is zero, leading to
logarithmic divergence.
We have adopted two strategies to address this difficulty: (i) removal of duplicates
after bootstrap resampling but prior to evaluating the estimator, and (ii) constructing an
ensemble of nearest-neighbour distances from the entire sample and bootstrapping on the
nearest-neighbour ensemble.
While these mitigations have different draw-backs, they are observed to perform equally well
when tested on a large variety of simple distributions and compared with accurate results for
the confidence intervals obtained from repeated Monte Carlo simulation trials.

Figure~\ref{BPEntropySim} presents the BP differential entropy versus traffic intensity $r$
for various server numbers $N$,
with overlaid simulation data as indicated by the dots with error bars.
The green dashed curve is the extreme large-$N$ limit based on a single exponential given by
\mbox{$\bar{F}(t) \asym{N\gg 1} e^{-t/m_{\text{bp}}}$}
with the mean $m_{\text{bp}}$ computed for $N = 60$.
The dashed black curve is the finer large-$N$ asymptotic form based on the two-exponential mixture
given in (\ref{SFAsymN}) and computed for $N = 60$.
It aligns very closely with the exactly calculated result.
Figure~\ref{BPLogMeanSim} presents the BP log-mean versus traffic intensity $r$
for various server numbers $N$,
with overlaid simulation data as indicated by the dots with error bars.
The green dashed curve is the extreme large-$N$ limit based on a single exponential given by
\mbox{$\bar{F}(t) \asym{N\gg 1} e^{-t/m_{\text{bp}}}$}
with the mean $m_{\text{bp}}$ computed for $N = 50$.
The dashed black curve is the finer large-$N$ asymptotic form based on the two-exponential mixture
given in (\ref{SFAsymN}) and computed for $N = 50$.
It aligns very closely with the exactly calculated result.

\section{Conclusions}
\label{Conclusions}
This paper has developed two distinct methods for generating explicit exact results for the
distribution of the partial BP pertaining to the M/M/$c$ queue and a variety of models
that generalize it with priority levels and distinct arrival classes.
The spectral method allows for a robust and efficient numerical implementation.
The algebraic method furnishes closed-form results for the PDF and SF that
elucidate the analytical structure of the problem.
In particular, for any given number of servers, it has identified a unique polynomial
that completely characterizes the distribution.
The present discussion has also served to connect previous diverse approaches
to the problem.

%
%
%

\begin{APPENDIX}{MGF Recurrence}
Let the RV
\mbox{$\mathcal{N}_{\text{b}}(t) = 0,1,\ldots,N$}
denote the number of servers that are busy at time $t$.
Let the RVs $T_n$ ,
\mbox{$n = 1,2,\ldots,N$},
denote the unit descent times
\begin{equation}
T_n \equiv \min\left\{t: \mathcal{N}_{\text{b}}(t) = k-1 | \mathcal{N}_{\text{b}}(0^+) = k\right\} \;,
\end{equation}
with their MGFs denoted by
\mbox{$\eta_n(s) \equiv \langle e^{-sT_n}\rangle_{T_n}$}.
\citet{BP:Omahen78} have shown that, for
\mbox{$n = 1,2,\ldots,N-1$},
\begin{equation}
\eta_n(s) = \frac{n\mu}{s + n\mu + \lambda -\lambda\eta_{n+1}(s)} \;, \quad
\eta_N(s) = \frac{N\mu}{s + N\mu + \lambda -\lambda\eta_N(s)} \;.
\end{equation}
This may be rearranged as
\begin{equation}
\mu_n - (s + \lambda + \mu_n)\eta_n(s) + \lambda\eta_n(s)\eta_{n+1}(s) = 0 \;,
\end{equation}
where
\mbox{$\mu_n \equiv n\mu$}.
Upon multiplying both sides by the product
\mbox{$\eta_1(s)\cdots\eta_{n-1}(s)$},
and introducing the RV convolution
\begin{equation}
\phi_n(s) = \prod_{k = 1}^n \eta_k(s) \;,
\label{EtaProd}
\end{equation}
we obtain
\begin{equation}
\mu_n\phi_{n-1}(s) - (s + \lambda + \mu_n)\phi_{n}(s) + \lambda\phi_{n+1}(s) = 0 \;,
\label{PhiRecur2}
\end{equation}
consistent with (\ref{PhiRecur}).
Also, we have from (\ref{EtaProd}) that
\mbox{$\phi_N(s) = \phi_{N-1}(s)\eta_N(s)$},
noting that
\mbox{$\eta_N(s)$}
is a solution of the quadratic equation
\begin{equation}
\lambda\eta^2_N(s) - (s + \mu_N + \lambda)\eta_N(s) + \mu_N = 0 \;,
\end{equation}
and is given by
\mbox{$\eta_N(s) = \psi_-(s)$},
with the minus branch chosen because
\mbox{$\psi_-(0) = 1$},
whereas
\mbox{$\psi_+(0) = 1/r > 1$}.
This completes the connection with (\ref{PhiRecur}).
By construction,
\mbox{$\eta_1(s)$}
is the MGF of the partial BP.
From the identification
\mbox{$\phi_1(s) = \eta_1(s)$},
it follows that, independent of the interpretation of the quantities
\mbox{$\phi_k(s)$},
\mbox{$k = 2,3,\ldots$},
the quantity
\mbox{$\phi_1(s)$},
derived as the solution of the system (\ref{PhiRecur2}),
also represents the MGF of the partial BP.
\end{APPENDIX}


%
\section*{Acknowledgments.}
The authors gratefully acknowledge useful discussions with Dr.~Stephen Bocquet.


\bibliographystyle{informs2014} 
\bibliography{BusyPeriod} 



\end{document}